\documentclass[12pt,twoside]{article}
\usepackage{amsmath}
\usepackage{amssymb}
\usepackage{amsthm}
\usepackage[latin1]{inputenc}
\usepackage[T1]{fontenc}
\usepackage{tabularx}
\usepackage{oldgerm}
\usepackage[all,2cell,ps]{xy}
\UseAllTwocells
\UsePSspecials{}
\usepackage{a4}

\newtheorem{thm}{Theorem}[subsection]
\newtheorem{prop}[thm]{Proposition}
\newtheorem{lemma}[thm]{Lemma}
\newtheorem{cor}[thm]{Corollary}

\theoremstyle{definition}

\newtheorem{nota}[thm]{Notation}

\newtheorem{defi}[thm]{Definition}
\newtheorem{remark}[thm]{Remark}

\def\df[#1]{%
  \ar@{}[#1]|(.35)*{}="A" \ar@{}[#1]|(.65)*{}="B" 
  \ar@{=>}"A";"B" }
\def\tdf[#1]{%
  \ar@{}[#1]|(.25)*{}="A" \ar@{}[#1]|(.55)*{}="B" 
  \ar@{=>}"A";"B" }

\newcommand{\on}{\operatorname}

\newcommand{\coker}{\operatorname{coker}}

\newcommand{\op}{\operatorname{op}}
\newcommand{\pt}{\operatorname{pt}}
\newcommand{\ol}{\overline}

\newcommand{\wtilde}{\widetilde}

\newcommand{\geqs}{\geqslant}
\newcommand{\leqs}{\leqslant}

\newcommand{\me}{\matheus}

\newcommand{\ra}{\rightarrow}
\newcommand{\lra}{\longrightarrow}
\newcommand{\isoto}{\overset{\sim}{\ra}}
\renewcommand{\to}[1][]{\xrightarrow{#1}}
\renewcommand{\isoto}[1][]{\xrightarrow[#1]{\sim}}

\newcommand{\Z}{\mathbb{Z}}

\newcommand{\ad}{\operatorname{ad}}
\newcommand{\Pic}{\operatorname{Pic}}

\DeclareMathAlphabet{\matheus}{U}{eus}{m}{n}
\newdir{ (}{{}*!/-5pt/@^{(}}
\newdir{ )}{{}*!/-5pt/@_{(}}
\numberwithin{equation}{subsection}

\newcommand{\id}{\operatorname{id}} 
\newcommand{\Ob}{\operatorname{Ob}} 
\newcommand{\Hom}[1][]{\on{Hom}_{\raise1.5ex\hbox to.1em{}#1}}
\newcommand{\Aut}[1][]{\operatorname{Aut}_{\raise1.5ex\hbox to.1em{}#1}} 

\newcommand{\cat}{\mathsf}

\newcommand{\catc}{{\mathsf{C}}} 
\newcommand{\catd}{{\mathsf{D}}} 
\newcommand{\catg}{{\mathsf{G}}} 
\newcommand{\cati}{{\mathsf{I}}} 

\newcommand{\catHom}[1][]{\mathsf{Hom}_{\raise1.5ex\hbox to.1em{}#1}}
\newcommand{\catEnd}[1][]{\mathsf{End}_{\raise1.5ex\hbox to.1em{}#1}}
\newcommand{\catAut}[1][]{\mathsf{Aut}_{\raise1.5ex\hbox to.1em{}#1}}
\newcommand{\catMod}{\mathsf{Mod}} 
\newcommand{\catTors}{\mathsf{Tors}} 
\newcommand{\catBitors}{\mathsf{Bitors}} 
\newcommand{\catSet}{\mathsf{Set}} 
\newcommand{\catGamma}{\mathsf{\Gamma}}
\newcommand{\catF}{\mathsf{F}}
\newcommand{\catQ}{\mathsf{Q}}

\newcommand{\catPSh}{\mathsf{PSh}} 
\newcommand{\catSh}{\mathsf{Sh}} 
\newcommand{\catCSh}{\mathsf{CSh}} 
\newcommand{\catOp}{\mathsf{Op}} 
\newcommand{\catGr}{\mathsf{Gr}} 
\newcommand{\catGrCom}{\mathsf{Gr^c}} 
\newcommand{\catMon}{\mathsf{Mon}} 
\newcommand{\catMonCom}{\mathsf{Mon^c}} 
\newcommand{\twocolim}[2]{\underset{#1}{2\!\varinjlim}{#2}} 
\newcommand{\threecolim}[2]{\underset{#1}{3\!\varinjlim}{#2}} 
\newcommand{\twolim}[2]{\underset{#1}{2\!\varprojlim}{#2}} 

\newcommand{\pstks}{\mathfrak{S}} 
\newcommand{\pstkt}{\mathfrak{T}} 

\newcommand{\stkg}{\mathfrak{G}} 
\newcommand{\stks}{\mathfrak{S}} 
\newcommand{\stkt}{\mathfrak{T}} 
\newcommand{\stkd}{\mathfrak{D}} 

\newcommand{\stkSh}{\mathfrak{Sh}} 
\newcommand{\stkLoc}{\mathfrak{LcSh}} 
\newcommand{\stkHom}[1][]{\mathfrak{Hom}_{\raise1.5ex\hbox to.1em{}#1}}
\newcommand{\stkEqui}[1][]{\mathfrak{Equi}_{\raise1.5ex\hbox to.1em{}#1}}
\newcommand{\stkTors}{\mathfrak{Tors}} 
\newcommand{\stkBitors}{\mathfrak{Bitors}} 
\newcommand{\stkAut}[1][]{\mathfrak{Aut}_{\raise1.5ex\hbox to.1em{}#1}}

\newcommand{\twostks}{\pmb{\mathfrak{S}}} 

\newcommand{\twostkSt}{\pmb{\mathfrak{St}}} 
\newcommand{\twostkLoc}{\pmb{\mathfrak{LcSh}}} 

\newcommand{\Cat}{\mathbf}
\newcommand{\Catc}{\mathbf{C}}
\newcommand{\CatSt}{\mathbf{St}} 
\newcommand{\CatPSt}{\mathbf{PSt}} 
\newcommand{\CatLcSt}{\mathbf{LcSt}} 
\newcommand{\CatCSt}{\mathbf{CSt}} 
\newcommand{\CatCat}{\mathbf{Cat}} 
\newcommand{\CatGr}{\mathbf{Gr}} 
\newcommand{\CatMon}{\mathbf{Mon}} 
\newcommand{\CatHom}[1][]{{\mathbf{Hom}}_{\raise1.5ex\hbox to.1em{}#1}}
\newcommand{\CatTop}{\mathbf{Top}} 
\newcommand{\CatGamma}{\mathbf{\Gamma}}

\newcommand{\TwoCat}{2\CatCat}
\newcommand{\TwoGr}{2\CatGr}
\newcommand{\TwoTop}{2\CatTop}
\newcommand{\TwoSt}{2\CatSt}
\newcommand{\TwoPSt}{2\CatPSt}
\newcommand{\Two}{\pmb}

\newcommand{\shhom}[1][]{{\me{H}om}_{\raise1.5ex\hbox to.1em{}#1}} 
\newcommand{\shaut}[1][]{{\me{A}ut}_{\raise1.5ex\hbox to.1em{}#1}} 

\newcommand{\pshf}{\me{F}} 
\newcommand{\pshg}{\me{G}} 

\newcommand{\she}{\me{E}} 
\newcommand{\shf}{\me{F}} 
\newcommand{\shg}{\me{G}} 
\newcommand{\shn}{\me{N}} 
\newcommand{\shp}{\me{P}} 
\newcommand{\shv}{\me{V}} 

\newcommand{\lshf}{\me{F}} 
\newcommand{\lshg}{\me{G}} 

\begin{document}

\title{Higher monodromy}

\author{P. Polesello and I. Waschkies} 

\date{}

\maketitle

\begin{abstract}
For a given category $\catc$ and a topological space $X$, the constant stack 
on $X$ with stalk $\catc$ is the stack of locally constant sheaves with values 
in $\catc$. Its global objects are classified by their monodromy, a functor 
from the Poincar\'e groupoid $\Pi_1(X)$ to $\catc$. In this paper we recall 
these notions from the point of view of higher category theory and then define 
the 2-monodromy of a locally constant stack with values in a 2-category 
$\Catc$ as a 2-functor from the homotopy 2-groupoid $\Pi_2(X)$ to $\Catc$. 
We show that 2-monodromy classifies locally constant stacks on a reasonably 
well-behaved space $X$. As an application, we show how to recover from this 
classification the cohomological version of a classical theorem of Hopf, 
and we extend it to the non abelian case.\\

\end{abstract}

\tableofcontents

\section*{Introduction}

A classical result in algebraic topology is the classification of the 
coverings of a (reasonably well-behaved) path-connected topological space $X$ 
by means of representations of its fundamental group $\pi_1(X)$.
In the language of sheaves this reads as an equivalence between the 
category of locally constant sheaves of sets on $X$ and that of 
representations of $\pi_1(X)$ on the stalk. The equivalence is given by the 
functor which assignes to each locally constant sheaf $\shf$ with stalk $S$ 
its monodromy $\mu(\shf)\colon\pi_1(X)\to\on{Aut}(S)$.

Now, let $\catc$ be a category. It makes sense to consider the monodromy
representations in $\catc$, in other words functors from the 
Poincar\'e groupoid $\Pi_1(X)$ to $\catc$. One would then say that they 
classify ``locally constant sheaves on $X$ with stalk in the category 
$\catc$'' even if there do not exist any sheaves with values in $\catc$. 
To state this assertion more precisely, one needs the language of stacks of 
Grothendieck and Giraud. A stack is, roughly speaking, a sheaf of categories 
and one may consider the constant stack $\catc_X$ on $X$ with 
stalk the category $\catc$ (if $\catc$ is the category of sets,
one recovers the stack of locally constant sheaves of sets). 
Then one defines a local system on $X$ with values in $\catc$ to be 
a global section of the constant stack $\catc_X$. The monodromy functor
establishes (for a locally relatively 1-connected space $X$) an equivalence of 
categories between global sections of $\catc_X$ and functors 
$\Pi_1(X)\lra\catc$.

A question naturally arises: what classifies stacks on $X$ which are 
locally constant? Or, on the other side, which geometrical objects are 
classified by representations (\emph{i.e.} 2-functors) of the homotopy 
2-groupoid $\Pi_2(X)$? In this paper, we define a locally constant stack 
with values in a 2-category $\Catc$ as a global section of the constant 
2-stack $\Catc_X$ and give an explicit construction of the 2-monodromy of 
such a stack as a 2-functor $\Pi_2(X) \lra \Catc$. We will show that, 
for locally relatively 2-connected topological spaces, a locally constant 
stack is uniquely determined (up to equivalence) by its 2-monodromy. 
We then use this result to recover the cohomological version of a classical 
theorem of Hopf, relating the second cohomology group with constant 
coeffiecients of $X$ to its first and second homotopy group, and we extend it 
to the non abelian case.

During the preparation of this work, a paper \cite{Toen} of B. Toen appeared,
where a similar result about locally constant $\infty$-stacks and 
their $\infty$-monodromy is established. His approach is different from ours, 
since we do not use any model category theory and any simplicial techniques, 
but only classical 2-category (and enriched higher category) theory. Moreover, 
since we are only interested in the degree 2 monodromy, we need weaker 
hypothesis on the space $X$ than \emph{loc.cit.}, where the author works on 
the category of pointed and connected $CW$-complexes.  

This paper is organised as follows. In Chapter \ref{locsheaves} we 
recall some basic notions of stack theory and give a functorial construction 
of the classical monodromy. Our approach appears at first view to be rather 
heavy, as we use more language and machinery in our definition as is usually 
done when one considers just monodromy for sheaves of sets  or abelian groups.
The reason for our category theoretical approach is to motivate the
construction of 2-monodromy (and to give a good idea how one could define
$n$-monodromy of a locally constant $n$-stack with values in an $n$-category,
for all $n$). As a byproduct we get the classification of locally constant 
sheaves with values in finite categories (\emph{e.g.} in the category defined 
by a group) which yields an amusing way to recover some non abelian 
versions of the "Hurewicz's formula", relating the first non abelian 
cohomology set with constant coefficients to representations of the 
fundamental group.

In Chapter \ref{locstacks} we introduce the 2-monodromy 2-functor of a locally 
constant stack with values in a 2-category. This construction is analogous to 
our approach to 1-monodromy, but the diagrams which should be checked for 
commutativity become rather large.
One reason for our lengthy tale on 1-monodromy is to give good evidence to
believe in our formulae, since we do not have the space to write down detailed 
proofs. We also describe the 2-monodromy as a descent datum on the loop space 
at a fixed point. Finally, we give some explicit calculations about the 
classification of gerbes with locally constant bands. This is related to 
Giraud's second non abelian cohomology set with constant coeffiecients. 

In Appendix \ref{sheavesC} we review the definition of the stack of 
sheaves with values in a complete category and in Appendix \ref{stacksC} the 
definition of the 2-stack of stacks with values in a 2-complete 2-category.

\noindent
{\bf Aknowledgement.}
We wish to thank Denis-Charles Cisinski for useful discussions and insights.

\section*{Notations and conventions}\label{notations}

We assume that the reader is familiar with the basic notions of
classical category theory, as those of category\footnote{There are a few 
well-known set-theoretical problems that arise in the deinition of a category. 
A convenient way to overcome these difficulties has been proposed by 
Grothendieck using the notion of a universe, and without saying so explicitly, 
we will work in this framework.}, functor between categories, 
transformation between functors (also called morphism of functors),
equivalence of categories, monoidal category and monoidal functor. We will 
also use some notions from higher category theory, as 2-categories, 2-functors,
2-transformations, modifications, 2-limits and 2-colimits. Moreover, we will 
look at 3-categories, 3-functors, etc., but only in the context of a
``category enriched in 2-categories'' which is much more elementary than the 
general theory of $n$-categories for $n\geqs 3$. 
References are made to \cite{ Borceux,Leinster,MacLane}\footnote{Note 
that 2-categories are called bicategories by some authors, for whom a 
2-category is what we call here a strict 2-category. A 2-functor is sometimes 
called a pseudo-functor.}. For an elementary introduction to 2-limits and 
2-colimits, see \cite{Waschkies}.

We use the symbols $\catc,\catd$, etc., to denote categories. 
If $\catc$ is a category, we denote by $\Ob \catc$ (resp. $\pi_0(\catc)$) the 
collection of its objects (resp. of isomorphism classes of its objects), and 
by $\Hom[\catc](P,Q)$ the set of morphisms between the objects $P$ and $Q$ 
(if $\catc=\catSet$, the category of all small sets, we will write
$\Hom(P,Q)$ instead of $\Hom[\catSet](P,Q)$). 
For a category $\catc$, its opposite category is denoted by $\catc^{\op}$.

We use the symbols $\Catc,\Cat{D}$, etc., to denote 2-categories. 
If $\Catc$ is a 2-category, we denote by $\Ob \Catc$ (resp. 
$\Two{\pi}_0(\Catc)$) the collection of its objects (resp. of equivalence 
classes of its objects), and by $\catHom[\Catc](\on{P},\on{Q})$ the small 
category of 1-arrows between the objects $\on{P}$ and $\on{Q}$ (if 
$\Catc=\CatCat$, the strict 2-category of all small categories, we 
will use the shorter notation $\catHom(\catc,\catd)$ to denote the category of 
functors between $\catc$ and $\catd$). Given two 2-categories $\Catc$ and 
$\Cat{D}$, the 2-category of 2-functors from $\Catc$ to $\Cat{D}$ will be 
denoted by $\CatHom (\Catc,\Cat{D})$. If $G$ is a commutative group, we will 
use the notation $\CatCat_G$ for the 2-category of $G$-linear categories and 
$\catHom[G](\catc,\catd)$ for the $G$-linear category of $G$-linear functors.  
If $\Catc$ is a 2-category, $\Catc^{\op}$ denotes its opposite 2-category, 
that is $\catHom[\Catc^{\op}](\on{P},\on{Q})=\catHom[\Catc](\on{Q},\on{P})$. 
Moreover, for any object $\on{Q}\in\Ob\Catc$, we set\footnote{This is 
consistent with the classical notion of Picard group. Indeed, if $R$ is a ring 
and $\catMod (R)$ the category of its left modules, by the Morita theorem the 
group $\Pic_{\CatCat_\Z}(\catMod(R))$ is isomorphic to the 
Picard group of $R$.} $\Pic_{\Catc}(\on{Q})=\pi_0(\catAut[\Catc](\on{Q}))$, 
where $\catAut[\Catc](\on{Q})$ denotes the monoidal category of 
auto-equivalences (\emph{i.e} 1-arrows which are invertible up to a 2-arrow) 
of $\on{Q}$ in $\Catc$. Note that the group $\Pic_{\Catc}(\on{Q})$ acts on the 
commutative group $\on{Z}_{\Catc}(\on{Q})=\
\Aut[{\catAut[\Catc]}(\on{Q})](\id_{\on{Q}})$ by conjugation.
If there is no risk of confusion, for a category (resp. $G$-linear category) 
$\catc$, we will use the shorter notations $\Pic (\catc)$ (resp. 
$\Pic_G(\catc)$) and $\on{Z}(\catc)$. This last group is usually called the 
center of $\catc$.

\section{Locally constant sheaves with values in a category}\label{locsheaves}

\subsection{Locally constant sheaves and their operations}

Let $X$ be a topological space. We denote by $\CatPSt(X)$ (resp. $\CatSt(X)$) 
the strict 2-category of prestacks (resp. stacks) of categories on $X$.
For the definition of a stack on a topological space, see 
Appendix \ref{stacksC}. 
The classical references are Giraud's book \cite{Giraud} and 
\cite[expos\'e VI]{SGA1}, and a more recent one is \cite{Brylinski}. 
A lighter presentation may be also found in \cite{D'Agnolo-Polesello}.

Recall that, as for presheaves, to any prestack $\stks$ on $X$ one naturally 
associates a stack $\stks^{\ddagger}$. Precisely, one has the following

\begin{prop}
     The forgetful 2-functor
     $$ \cat{For}:\ \CatSt(X) \lra \CatPSt(X) $$
     has a 2-left adjoint 2-functor
     $$ \ddagger:\ \CatPSt(X) \lra \CatSt(X). $$
\end{prop}
\noindent 
Let us fix an adjunction 2-transformation
$$ \eta_X:\ \cat{Id}_{\CatPSt(X)} \lra \cat{For}\circ \ \ddagger. $$ 

Note that there is an obvious fully faithful\footnote{Recall that a 2-functor 
$\catF\colon\Catc\to\Cat{D}$ is faithful (resp.\ full, resp.\ fully faithful)
if for any objects $\on{P},\on{Q}\in \Ob\Catc$, the induced functor 
$\catF\colon \catHom[\Catc](\on{P},\on{Q})\to \catHom[\Cat{D}](\catF(\on{P}),
\catF(\on{Q}))$ is faithful (resp.\ full and essentially surjective, 
resp.\ an equivalence of categories).}
2-functor of 2-categories $\CatCat \lra \CatPSt(X)$ which associates to a 
category $\catc$ the constant prestack on $X$ with stalk $\catc$.

\begin{defi}
     Let $\catc$ be a category. The constant stack on $X$ with stalk $\catc$ 
     is the image of $\catc$ by the 2-functor
     $$ 
     (\,\cdot\,)_X\colon \CatCat \lra \CatPSt(X) \overset{\ddagger}{\lra} 
     \CatSt(X).
     $$   
\end{defi}
\noindent
Note that the 2-functor $ (\,\cdot\,)_X$ preserves faithful and fully faithful 
functors (hence sends (full) subcategories to (full) substacks). Moreover, the 
2-transformation $\eta_X$ induces on global sections a natural faithful functor
$$ \eta_{X,\catc}\colon\catc \lra \catc_X(X). $$

\begin{defi}
     An object $\lshf\in\Ob{\catc_X}(X)$ is called a local system on 
     $X$ with values in $\catc$. A local system is constant with 
     stalk $M$ if it is isomorphic to $\eta_{X,\catc}(M)$ for some object 
     $M\in\Ob{\catc}$.
\end{defi}

Let $\catc=\catSet$, the category of all small sets (resp. $\catc=\catMod(A)$,
the category of left $A$-modules for some ring $A$). Then it is easy to see
that $\catc_X$ is naturally equivalent to the stack of locally constant 
sheaves of sets (resp. $A_X$-module). Moreover the functor 
$\eta_{X,\catc}\colon\catc\ra\catc_X(X)$ is canonically equivalent to the 
functor which associates to a set $S$ (resp. an $A$-module $M$) the constant 
sheaf on $X$ with stalk $S$ (resp. $M$). More generally one can easily prove the 
following proposition:

\begin{prop}
     Let $\catc$ be a complete category\footnote{Recall that a complete 
     category is a category admitting all small limits.} and $X$ a locally 
     connected topological space. Denote by $\stkLoc_X(\catc)$ the full 
     substack of the stack $\stkSh_X(\catc)$ of sheaves with values in $\catc$,
     whose objects are locally constant. Then there is a natural equivalence 
     of stacks
     $$ \catc_X  \overset{\sim}{\lra} \stkLoc_X(\catc).$$ 
\end{prop}
\noindent
For a detailed construction of the stacks  $\stkSh_X(\catc)$ and 
$\stkLoc_X(\catc)$, see Appendix \ref{sheavesC}.

\begin{remark}
Suppose now that $\catc$ is a category that is not necessarily complete. 
The  category $\widehat{\catc}=\catHom(\catc^{\op},\catSet)$ is complete 
(and cocomplete) and the Yoneda embedding
$$ Y\colon\catc \lra \widehat{\catc} $$ 
commutes to small limits. 
Then one usually defines sheaves with values in $\catc$ as presheaves 
that are sheaves in the category $\catPSh_X(\widehat{\catc})$. Note that 
in general there does not exist a sheaf with values in $\catc$ (take for 
example $\catc$ the category of finite sets), but if $\catc\neq \varnothing$ 
then the category $\catc_X(X)$ of local systems with values in $\catc$ is 
always non-empty. Still, even in that case, we will sometimes refer to
a local system as a locally constant sheaf. More precisely we get a fully 
faithful functor of stacks
$$Y_X\colon\catc_X \lra (\widehat{\catc})_X \simeq\stkLoc_X(\widehat{\catc}).$$
Let $\me{F}$ be an object of $\stkLoc_X(\widehat{\catc})$ on some open set. 
Then $\me{F}$ is in the essential image of $Y_X$ if and only if all of its 
stalks are representable.
\end{remark}

Let $f\colon X\ra Y$ be a continuous map of topological spaces, $\stks$ a 
prestack on $X$ and $\stkd$ a prestack on $Y$.

\begin{nota}
\begin{itemize}

     \item[(i)] Denote by $f_*\stks$ the prestack on $Y$ such that, for any 
     open set $V\subset Y$, $f_*\stks(V)= \stks(f^{-1}(V))$. 
     If $\stks$ is a stack on $X$, then $f_*\stks$ is a stack on $Y$.
 
     \item[(ii)] Denote by $f^{-1}_{p}\stkd$ the prestack on $X$ such that, 
     for any open set $U\subset X$, $f^{-1}_p\stkd(U)= \twocolim{f(U)\subset V}
     {\stkd(V)}$. If $\stkd$ is a stack on $Y$, we set $ f^{-1}\stkd = 
     (f^{-1}_p\stkd)^{\ddagger}$.
\end{itemize}
\end{nota}
\noindent
Recall that the category $\twocolim{f(U)\subset V}{\stkd(V)}$ is described as 
follows:
\begin{align*}
\Ob \twocolim{f(U)\subset V}{\stkd(V)} & = 
\bigsqcup_{f(U)\subset V}\Ob \stkd(V),\\
\Hom[{\twocolim{f(U)\subset V}{\stkd(V)}}](G_{V},G_{V'}) & = 
\underset{f(U)\subset V''\subset V\cap V'}{\varinjlim}
\Hom[\stkd(V'')](G_{V}|_{V''},G_{V'}|_{V''}).
\end{align*}

\begin{prop}\label{fadj}
The 2-functors
$$ f_*\colon \CatSt(X) \lra \CatSt(Y) \qquad f^{-1} \colon \CatSt(Y) \lra 
     \CatSt(X)$$
are 2-adjoint, $f_*$ being the right 2-adjoint of $f^{-1}$.   
\end{prop}
\noindent
Moreover, if $g\colon Y \to Z$ is another continuous map, one has natural 
equivalences\footnote{For sake of simplicity, here and in the sequel we use 
the word ``equivalence'' for a 2-transformation which is invertible up to a 
modification. In the case of the inverse image, to be natural means that if we 
consider the continuous maps $h=f_3\circ f_2\circ f_1$, $g_1=f_2\circ f_1$ and 
$g_2=f_3\circ f_2$, then the two equivalences 
$h^{-1}\simeq f_1^{-1}g^{-1}_2\simeq f^{-1}_1f^{-1}_2 f^{-1}_3$ and 
$h^{-1}\simeq g^{-1}_1f^{-1}_3\simeq f^{-1}_1f^{-1}_2f^{-1}_3$ 
are naturally isomorphic by a modification, in the sense that, if we look
at the continuous map $k=f_4\circ f_3\circ f_2\circ f_1$, we get the obvious 
commutative diagram of modifications.} of 2-functors
$$
g_*\circ f_*\simeq (g\circ f)_*, \qquad f^{-1}\circ  g^{-1}\simeq 
(g\circ f)^{-1}.
$$

Since the following diagram commutes up to equivalence
\begin{equation}\label{diagram1}
     \xymatrix@R=.3cm{ 
     &  {\CatPSt(Y)} \ar[dd]^{f^{-1}_p} \ar[r]^{\ddagger} &{\CatSt(Y)} 
     \ar[dd]^{f^{-1}}  \\
     {\CatCat} \ar[ur] \ar[dr]   & &  \\
     &  {\CatPSt(X)} \ar[r]^{\ddagger} & {\CatSt(X)}  ,} 
\end{equation}
the 2-functor $f^{-1}$ preserves constant stacks (up to natural equivalence).

\begin{defi}
     Denote by $\catGamma(X,\,\cdot\,)$ the 2-functor of global sections
     $$\CatSt(X)\lra\CatCat \quad ; \quad \stks\mapsto 
     \catGamma(X,\stks)=\stks(X) $$
     and set $\catGamma_X=\catGamma(X,\cdot)\circ (\cdot)_X$. 
\end{defi}

Note that for any stack $\stks$, the category $\stks(\varnothing)$ is the 
terminal object in $\CatCat$ (which consists of precisely one morphism). 
Hence the 2-functor
$$\catGamma (\{\pt\},\cdot):\ \CatSt(\{\pt\})\lra\CatCat $$
is an equivalence of 2-categories. We easily deduce

\begin{prop}
     The 2-functor $\catGamma(X,\,\cdot\,)$ is right 2-adjoint to 
     $(\,\cdot\,)_X$.
\end{prop}

\noindent
It is not hard to see that we may choose the functors 
$\eta_{X,\catc}\colon\catc \lra \catGamma(X,\catc_X)$ to 
define the adjunction 2-transformation
$$\eta_X \colon \cat{Id}_\CatCat \lra \catGamma_X.$$

Consider the commutative diagram of topological spaces
$$ 
\xymatrix@R=.4cm@C=.5cm{
       X \ar[rr]^f \ar[dr]_{a_X}  & & Y \ar[dl]^{a_Y}  \\
       & \{\pt\} & } 
$$
and the induced 2-transformation of 2-functors
$$ 
a_{Y*}\circ (\cdot)_Y \lra a_{Y*}\circ f_* \circ f^{-1}\circ (\cdot)_Y 
\simeq a_{X*}\circ f^{-1}\circ (\cdot)_Y \simeq a_{X*}\circ (\cdot)_X .
$$
Hence we get a 2-transformation of 2-functors $f^{-1}$ compatible with 
$\eta_X$ and $\eta_Y$, \emph{i.e.} such that the following diagram commutes up
to a natural invertible modification:
$$ 
\xymatrix@R=.4cm@C=.5cm{
       {\catGamma_Y} \ar[rr]^{f^{-1}} & & {\catGamma_X} .\\
       & {\cat{Id}_\CatCat }\ar[ul]^{\eta_Y} \ar[ur]_{\eta_X} & } 
$$
Note that this implies that for each point $x\in X$ and any local systems 
$\shf,\shg\in \catc_X(X)$, the natural morphism
$$ 
i_x^{-1}\shhom[\catc_{X}](\shf,\shg) \ra \Hom[\catc](i^{-1}_x\shf,i_x^{-1}
\shg) 
$$ 
is an isomorphism (here $i_x\colon \{\pt\}\to X$ denotes the natural map 
sending $\{\pt\}$ to $x$ and we identify $\catc$ with global sections of 
$\catc_{\pt}$). Therefore, for any continuous map $f\colon X\ra Y$, we get a 
natural isomorphism of locally constant sheaves of sets
$$ 
f^{-1}\shhom[\catc_{Y}](\shf,\shg)  \isoto
 \shhom[\catc_X](f^{-1}\shf,f^{-1}\shg). 
$$ 

Now we are ready to formulate the fundamental Lemma of homotopy invariance
of local systems as follows:

\begin{lemma}\label{XxI}
     Let $I=[0,1]$ and $t\in I$. Consider the maps
     $$
     \xymatrix@C=2em{{X} \ar@<.5ex>[r]^-{\iota_t} & 
     {X\times I} \ar@<.5ex>[l]^-{p} ,}
     $$
     where $\iota_t(x)=(x,t)$ and $p$ is the projection. Then the 
     2-transformations 
     $$
     \xymatrix@C=3em{ {\catGamma_X} \ar@<.5ex>[r]^-{p^{-1}} 
     &{\catGamma_{X\times I}} \ar@<.5ex>[l]^-{\iota_t^{-1}} }
     $$
     are equivalences of 2-functors, quasi-inverse one to each other. 
\end{lemma}

\begin{proof}
Let $\catc$ be a category. It is sufficient to show that the functors
$$
\xymatrix@C=3em{ {\catGamma(X,\catc_X)} \ar@<.5ex>[r]^-{p^{-1}} 
&{\catGamma(X\times I,\catc_{X\times I})} \ar@<.5ex>[l]^-{\iota_t^{-1}}}
$$
are natural equivalences of categories, quasi-inverse one to each other.
Since $\iota_t^{-1}\circ p^{-1}\simeq (p\circ \iota_t)^{-1} = 
\id_{\catGamma(X,\catc_X)}$, it remains to check that for each 
$\shf\in\catGamma(X\times I,\catc_{X\times I})$ there is a natural isomorphism 
$p^{-1} \iota_t^{-1}\shf\simeq \shf$.\\
First let us prove that if $\shf$ is a locally constant sheaf of sets on 
$X\times I$, then the natural morphism
\begin{equation}\label{XxI,F}
\iota_t^{-1}\colon \Gamma(X\times I,\shf) \ra \Gamma(X,\iota_t^{-1}\shf)
\end{equation}
is an isomorphism. Indeed, let $s$ and $s'$ be two sections in 
$\Gamma(X\times I,\shf)$ such that $s_{x,t} = s'_{x,t}$. 
Since $\shf$ is locally constant, the set $\{t'\in I\ |\ s_{x,t'} = s_{x,t}\}$
is open and closed, hence equal to $I$. Therefore the map is injective. 
Now let $s\in \Gamma(X,i_t^{-1}\shf)$. Then $s$ is given by sections $s_j$ of 
$\shf$ on on a family $(U_j\times I_j)_{j\in J}$ where $I_j$ is an open 
interval containing $t$, the $U_j$ cover $X$ and the sheaf $\shf$ is constant 
on $U_j\times I_j$. It is not hard to see that, by refining the covering,
the sections $s_j$ can be extended to $U_j\times I$ and using the injectivity 
of the map one sees that we can patch the extensions of the $s_j$ to get a 
section of $\shf$ on $X\times I$ that is mapped to $s$. Hence the morphism 
\eqref{XxI,F} is an isomorphism.\\
Now let $\shf\in\catGamma(X\times I,\catc_{X\times I})$.
Since $\iota^{-1}_t \shhom[\catc_{X\times I}](p^{-1}\iota_{t}^{-1}\shf,\shf)
\isoto \shhom[\catc_{X}](\iota^{-1}_t\shf,\iota^{-1}_t\shf)$, we get the 
isomorphism
$$ 
\Gamma(X\times I,\shhom[\catc_{X\times I}](p^{-1}\iota_{t}^{-1}\shf,\shf))
\isoto \Gamma(X,\shhom[\catc_{X}](\iota^{-1}_t\shf,\iota^{-1}_t\shf)).
$$ 
It is an easy exercise to verify that we can defines the isomorphism 
$p^{-1}\iota_{t}^{-1}\shf \isoto \shf$ by the identity section in the set 
$\Gamma(X,\shhom[\catc_X](\iota_t^{-1}\shf,\iota_{t}^{-1}\shf))$. 
\end{proof}

\noindent
Taking $X=\{\pt\}$ in Lemma \ref{XxI}, we get

\begin{cor}
     The adjunction 2-transformation 
     $$ \eta_I\colon \cat{Id}_{\CatCat} \lra \catGamma_I $$
     is an equivalence, i.e. for each category $\catc$ the functor
     $\eta_{I,\catc}\colon \catc \lra \catGamma(I,\catc_I)$ is a natural 
     equivalence.
\end{cor}

\noindent
For each $X$ and each $t\in I$, we have an invertible modification 
$\id_{\catGamma_X}\simeq \iota_t^{-1}p^{-1}$. Therefore, for any $s,t\in I$ 
there exists a canonical invertible modification $\iota_t^{-1}\simeq \iota_s^{-1}$
which can be used to prove the following technical Lemma:

\begin{lemma}\label{auxlemma1}
     The diagram of continuous maps on the left induces for any $s,t,t'\in I$ 
     the commutative diagram of modifications on the right
     $$ 
     \xymatrix{ X \ar[r]^{\iota_t} \ar[d]_{\iota_s} & X\times I
     \ar[d]^{\iota_s\times\id_I} \\
     X\times I \ar[r]_{\id_X\times\iota_t} & X\times I^2 } 
     \qquad \qquad
     \xymatrix{ \iota_t^{-1}(\iota_s\times\id_I)^{-1} \ar[r]^{\sim}
     \ar[d]_{\wr} & 
     \iota_{t'}^{-1}(\iota_s\times \id_I)^{-1}  \ar[d]^{\wr}  \\
     \iota_s^{-1}(\id_X\times\iota_t)^{-1} \ar[r]_{\sim} & 
     \iota_{s}^{-1}(\id_X\times\iota_{t'})^{-1} .}
     $$
     Moreover, let $H\colon X\times I\ra Y$ be a continuous map that factors 
     through the projection 
     $p\colon X\times I\ra X$. Then the composition of invertible modifications
     $$ 
     (H\circ \iota_{t})^{-1}\simeq \iota_{t}^{-1}H^{-1} \simeq 
     \iota_{t'}^{-1}H^{-1} \simeq (H\circ \iota_{t'})^{-1} 
     $$
     is the identity.
\end{lemma}

Let $\CatTop$ denote the strict 2-category of topological spaces and 
continuous maps, where 2-arrows are homotopy classes of homotopies between 
functions (see for example \cite[cap.~7]{Borceux} for explicit details). 
Then homotopy invariance of locally constant sheaves may be expressed 
as the following

\begin{prop}\label{the2functorgamma}
     The assignment $(\catc,X)\mapsto\catGamma(X,\catc_X)$ defines a 2-functor
     $$\catGamma\colon \CatCat\times\CatTop^{\op}\lra\CatCat. $$
     Moreover, the natural functors $\eta_{X,\catc}\colon 
     \catc\lra\catGamma(X,\catc)$ define a 2-transformation
     $$ \eta\colon\catQ_1 \lra \catGamma $$
     where $\catQ_1\colon \CatCat\times\CatTop^{\op} \lra \CatCat$ is the 
     projection.
\end{prop}

\begin{proof}
It remains to check that $\catGamma$ is well defined at the level of 2-arrows.
Let $f_i\colon X\to Y$ be continuous maps for $i=0,1$ and let $H\colon f_0\to 
f_1$ be a homotopy, that is a continuous map $X\times I\ra Y$ such 
that $H\circ\iota_0=f_0$ and $H \circ\iota_1=f_1$. 
Then $\alpha_H:f^{-1}_0\isoto f^{-1}_1$ is defined by the chain of natural 
invertible modifications
$$ 
f^{-1}_0=(H\circ \iota_0)^{-1}\simeq \iota_{0}^{-1}H^{-1} \simeq
\iota_1^{-1}H^{-1} \simeq (H\circ \iota_1)^{-1}=f_1^{-1}. 
$$
Consider the constant homotopy at $f\colon X\to Y$
$$ H_f\colon X\times I\ra Y \quad ; \quad (x,t)\mapsto f(x). $$ 
Since we may factor $H_f$ as $X\times I \to[p] X \to[f] Y$,
by Lemma \ref{auxlemma1} we get $ \alpha_{H_f}=\on{id_{f^{-1}}}$.\\
Now let $H_0,H_1\colon f_0\to f_1$ be two homotopies and 
$K\colon H_0\to H_1$ a homotopy. Then consider the commutative diagram
of invertible modifications
$$ 
\xymatrix{ f^{-1}_0 \ar[r]^-\sim \ar[d]_-{\alpha_{H_0}} & 
   \iota_{0}^{-1}H_0^{-1} \ar[r]^-\sim \ar[d]^-{\wr} & 
   \iota_0^{-1}j_{0}^{-1}K^{-1} \ar[r]^-{\sim} \ar[d]^-{\wr}  &
   \iota_0^{-1}j_{1}^{-1}K^{-1} \ar[d]^-{\wr} &  
   \iota_0^{-1}H_1^{-1} \ar[l]_-{\sim} \ar[d]^-{\wr} & 
   f_0^{-1} \ar[l]_-{\sim} \ar[d]^-{\alpha_{H_1}}  \\
   f^{-1}_1 \ar[r]_-{\sim} & \iota_{1}^{-1}H_0^{-1} \ar[r]_-{\sim} & 
   \iota_1^{-1}j_{0}^{-1}K^{-1} \ar[r]_-{\sim} &
   \iota_1^{-1}j_{1}^{-1}K^{-1}  &  \iota_1^{-1}H_1^{-1} \ar[l]^-{\sim} & 
   \ar[l]^-{\sim} f_1^{-1}.  }
$$ 
We have to check that the horizontal lines are identity modifications. 
This is a consequence of Lemma \ref{auxlemma1}, which allows us (by a diagram 
chase) to identify the two lines with the modifications induced by constant 
homotopies.\\
The fact that $\alpha$ is compatible with the composition of homotopies is 
finally a very easy diagram chase.
\end{proof}

\subsection{The monodromy functor}

\begin{defi}
     The homotopy groupoid (or Poincar\'e groupoid) of $X$ is the small 
     groupoid $\Pi_1(X)=\catHom[\CatTop](\{\pt\},X)$. 
\end{defi}

Roughly speaking, objects of $\Pi_1(X)$ are the points of $X$ and if 
$x,y\in X$, $\on{Hom}_{\Pi_1(X)}(x,y)$ is the set of homotopy classes of 
paths starting from $x$ and ending at $y$. The composition law is the 
opposite of the composition of paths. 
Note that, in particular, $\pi_0(\Pi_1(X)) = \pi_0(X)$, the set of arcwise 
connected components of $X$, and for each $x\in X$, $\Aut[\Pi_1(X)](x) = 
\pi_1(X,x)$, the fundamental group of $X$ at $x$. 

Let $\Cat{Gr}$ denote the strict 2-category of small groupoids. Then we 
have a 2-functor 
$$
\Pi_1\colon\CatTop\lra\Cat{Gr},
$$ 
which defines the Yoneda type 2-functor
$$ 
\cat{Y}_{\Pi_1}\colon \CatCat \times \CatTop^{\op} \lra \CatCat, \quad (\catc, X) \mapsto 
\catHom(\Pi_1(X),\catc) .
$$

\begin{defi}
     The monodromy 2-transformation
     $$ \mu\colon \catGamma \lra \cat{Y}_{\Pi_1} $$
     is defined as follows. For any topological space $X$ and any
     category $\catc$, the functor 
     $$ \catGamma_{X,\catc}\colon \catHom[\CatTop](\{\pt\},X) \lra
     \catHom(\catGamma(X,\catc_X),\catc) $$
     induces by evaluation a natural functor
     $$ \catGamma(X,\catc_X) \times\Pi_1(X) \lra \catc, $$
     hence by adjunction a functor
     $$ \mu_{X,\catc}\colon \catGamma(X,\catc_X)\lra \catHom(\Pi_1(X),\catc). $$
     We will sometimes use the shorter notation $\mu$ instead 
     of the more cumbersome $\mu_{X,\catc}$. We will also extend $\mu$ to 
     pointed spaces without changing the notations.
\end{defi}

Let us briefly illustrate that we have constructed the well-known classical 
monodromy functor.
Recall that, for each $x\in X$ one has a natural stalk 2-functor
$$\catF_x\colon \CatSt(X)\lra\CatCat \quad ; \quad \stks\mapsto \stks_x=
\twocolim{x\in U}{\ \stks(U).} $$ 
Let $i_x\colon \{x\}\to X$ denote the natural embedding. Since $\catF_x$ is 
canonically equivalent to $\catGamma(\{x\},\cdot)\circ i_x^{-1}$, one
gets a 2-transformation $\catGamma(X,\cdot)\to\catF_x$ and then a 
2-transformation 
$$\rho_x\colon \catGamma_X\lra \catF_{X,x}=\catF_x\circ (\cdot)_X .$$
For an object $\shf$ (resp. morphism $f$) in $\catc_X(X)$, if there is no 
risk of confusion, we will simply write $\shf_x$ (resp. $f_x$) to denote 
$\rho_{x,\catc}(\shf)$ (resp. $\rho_{x,\catc}(f)$) in $(\catc_X)_x$.

Let $\CatTop_*$ be the strict 2-category of pointed topological spaces, 
pointed continuous maps and homotopy classes of pointed homotopies. One can 
prove by diagram chases similar to those in the proof of 
Proposition~\ref{the2functorgamma}, that 

\begin{prop}
     The assignment $(\catc,(X,x)) \mapsto (\catc_X)_x$ defines a 2-functor
     $$ \catF\colon \CatCat \times \CatTop_*^{\op} \lra \CatCat .$$
     Moreover, the natural functors $\rho_{x,\catc}\colon \catGamma(X,\catc_X) 
     \lra (\catc_X)_x $ define a 2-transformation
     $$ \rho\colon \catGamma \lra \catF. $$
\end{prop}
Since the stalks of the stack associated to a prestack do not change, for 
each category $\catc$ and each pointed space $(X,x)$ the functor 
$$\catc\to[\eta_{X,\catc}] \catc_X(X) \to[\rho_{x,\catc}] (\catc_X)_x$$ 
is an equivalence of categories. Hence the composition
$$ \rho\circ\eta\colon \catQ_1 \lra \catF $$
is an equivalence of 2-functors (here $\catQ_1\colon\CatCat \times 
\CatTop_*^{\op} \lra \CatCat$ is the natural projection). 
Let $\varepsilon\colon\catF\lra\catQ_1$ denote a fixed quasi-inverse to 
$\rho\circ\eta$, \emph{i.e.} for each category $\catc$ and each pointed space 
$(X,x)$, we fix a natural equivalence 
$$ \varepsilon_{x,\catc}\colon (\catc_X)_x \overset{\sim}{\lra} \catc$$
such that if we have a pointed continuous map $f\colon(X,x)\ra (Y,y)$ we get a 
diagram
$$ \xymatrix{
     {\catc_Y(Y)} \ar[rr]^{f^{-1}} \ar[dd]_{\rho_{y, \catc}} & & 
     {\catc_X(X)} \ar[dd]^{\rho_{x, \catc}} \\
     & {\catc} \ar[ul]^{\eta_{Y, \catc}} \ar[ur]_{\eta_{X, \catc}}   & \\
     {(\catc_Y)_{y}} \ar[rr]^{\sim}_{f^{-1}} 
     \ar[ur]_{\varepsilon_{y,\catc}}^\sim   & & 
     {(\catc_X)_{x}} \ar[ul]^{\varepsilon_{x, \catc}}_{\sim} }$$
that commutes up to natural isomorphism.
Denote by $\omega$ the composition $\varepsilon \circ \rho \colon \catGamma 
\lra \cat{Q}_1$. Fix a topological space $X$ and a category $\catc$, and let 
$\shf\in\catGamma(X,\catc_X)$. Then a direct comparison shows that, up to 
natural isomorphism, we have
$$ \mu_{X,\catc}(\shf)(x)=\omega_{x,\catc}(\shf) $$
(if $\catc=\catSet$, then $\omega_x$ is just the usual stalk-functor)
and if $\gamma\colon x\ra y$ is a path, then $\mu_{X,\catc}(\shf)(\gamma)$ is 
defined by the chain of isomorphisms
$$ 
\omega_{x, \catc}(\lshf) \simeq \omega_{0, \catc}(\gamma^{-1}\lshf) \simeq 
\eta_{I, \catc}^{-1}(\gamma^{-1}\lshf) \simeq 
 \omega_{1, \catc}(\gamma^{-1}\lshf) \simeq \omega_{y, \catc}(\lshf) 
$$
(and if $\catc=\catSet$, we usually choose $\eta_I^{-1} = 
\Gamma(I,\,\cdot\,)$).\\
In particular, this means that the following diagram commutes up to natural 
invertible modification 
\begin{equation}\label{diagram2}\
\xymatrix@C=.7cm@R=.5cm{ {\catGamma}\ar[rr]^{\mu}  \ar[dr]_{\omega} && 
     {\cat{Y}_{\Pi_1},} \ar[dl]^{ev} \\
     & {\cat{Q}_1} & }
\end{equation}
where $ev$ is the evaluation 2-transformation, that is for each pointed space 
$(X,x)$ and each functor $\alpha\colon\Pi_1(X)\to\catc$, it is defined by 
$ev_x(\alpha) = \alpha(x)$.

Let $\Delta\colon\cat{Q}_1 \lra \cat{Y}_{\Pi_1}$ denote the diagonal 
2-transformation: for each topological space $X$, each category $\catc$ and 
each $M\in\Ob\catc$, $\Delta_{X,\catc} (M)$ is the constant functor 
$x\mapsto M$ (\emph{i.e.} the trivial representation with stalk $M$). 
Clearly $ev\circ\Delta = \id_{\cat{Q}_1}$. Moreover, by a diagram chase, we 
easily get

\begin{prop} 
     The diagram of 2-transformations
     $$ 
     \xymatrix@C=.7cm@R=.5cm{ {\catGamma}\ar[rr]^{\mu} & & {\cat{Y}_{\Pi_1}} \\
     & \ar[ul]^{\eta} {\cat{Q}_1} \ar[ur]_{\Delta} & }
     $$
     commutes up to invertible modification.  
\end{prop}

\begin{prop}\label{1-faith}
     For each topological space $X$ and each category $\catc$, the functor 
     $$ \mu_{X,\catc}\colon \catGamma(X,\catc_X) \lra\catHom(\Pi_1(X),\catc).$$
     is faithful and conservative.
\end{prop}

\begin{proof}
Let $f,g\colon\lshf\ra\lshg$ be two morphisms of local systems such 
that $\mu (f)=\mu (g)$. Since the diagram \eqref{diagram2} commutes and 
$\varepsilon_{x,\catc}\colon (\catc_X)_x \lra\catc$ is an equivalence, 
we get that $f_x=g_x$ in $(\catc_X)_x$ for all $x\in X$. 
Since $\catc_X$ is a stack, this implies that $f=g$, hence $\mu_{X, \catc}$ is 
faithful.\\
The same diagram implies that if $\mu (f)$ is an isomorphism, then
$f_x$ is an isomorphism in $(\catc_X)_x$ for all $x\in X$ and therefore $f$ is 
an isomorphism.
\end{proof}

\begin{prop}\label{1-full}
     Let $X$ be locally arcwise connected. Then for each category $\catc$, 
     the functor 
     $$\mu_{X,\catc}\colon \catGamma(X,\catc_X) \lra\catHom(\Pi_1(X),\catc).$$ 
     is full.
\end{prop}

\begin{proof}
Let $\lshf, \lshg\in\catc_X(X)$. A morphism $\phi\colon\mu(\lshf)\ra
\mu(\lshg)$ is given by a family of morphisms 
$$ \phi_x\colon ev_x(\mu(\lshf)) \to ev_x(\mu(\lshg)) $$
such that for every (homotopy class of) path $\gamma\colon x\ra y$, the diagram
$$ \xymatrix{
      ev_x(\mu(\lshf)) \ar[r]^{\phi_x} \ar[d]_{\mu(\lshf)(\gamma)} &  
      ev_x(\mu(\lshg)) \ar[d]^{\mu(\lshg)(\gamma)}  \\
      ev_y(\mu(\lshf)) \ar[r]_{\phi_y} & ev_y(\mu(\lshg)) } $$
is commutative. Using the diagram \eqref{diagram2} and the definition of the 
stalk of a stack, we get an arcwise connected open neighborhood $U_x$ of $x$ 
and a morphism 
$$ \varphi^x\colon \lshf|_{U_x} \ra \lshg|_{U_x} $$
such that $ev_x(\mu(\varphi^x))=\phi_x$.\\
In order to patch the $\varphi^x$, we have to show that for any 
$z\in U_x\cap U_y$ we have $(\varphi^x)_z=(\varphi^y)_z$. Since $\varepsilon_z$
is an equivalence, it is sufficent to check that 
$ev_z(\mu(\varphi^x))=ev_z(\mu(\varphi^y))$. \\
Choose a path $\gamma \colon x\rightarrow z$. Then
$$\mu(\lshg)(\gamma) \circ ev_x(\mu(\varphi^x)) \circ \mu(\lshf)(\gamma)^{-1}= 
   ev_z(\mu(\varphi^x)).  $$
But by definition the lefthand side is just $\phi_z$. Similarly, taking a path
$\gamma'\colon y \rightarrow z$, we get
$$  ev_z(\mu(\varphi^x))=\phi_z= ev_z(\mu(\varphi^y)) $$ 
and, by definition of the stalk, this means that $\varphi^x$ and $\varphi^y$ 
coincide in a neighborhood of $z$. Since $z$ was chosen arbitrarily, they 
coincide on $U_x\cap U_y$. Since $\catc_X$ is a stack, we can patch the 
morphisms $\varphi_x$ to a unique $\varphi\colon\lshf\ra\lshg$ such that 
$\mu(\varphi)=\phi$.
\end{proof}

\begin{cor}\label{1-con}
     Let $X$ be 1-connected \footnote{Here and for the sequel, a topological 
     space $X$ is $n$-connected if $\pi_i(X)\simeq 1$ for all 
     $0\leqs i\leqs n$, and locally $n$-connected if each neighborhood of each 
     point contains an $n$-connected neighborhood.}. Then 
     $$\eta_X \colon\cat{Id}_{\CatCat} \lra \catGamma_X$$ 
     is an equivalence of 2-functors, i.e. for each category $\catc$, 
     the functor $\eta_{X,\catc}\colon\catc \lra \catGamma(X, \catc_X)$ is an 
     equivalence of categories, natural in $\catc$.
\end{cor}

\begin{proof}
Fix $x_0\in X$. Since the groupoid $\Pi_1(X)$ is trivial (that is $\Pi_1(X)
\simeq 1$), the 2-transformation $ev_{x_0}\colon \cat{Y}_{\Pi_1(X)} \lra
\cat{Id}_{\CatCat}$ is an equivalence, quasi-inverse 
to $\mu_X\circ\eta_X$. Since $\mu_X$ is fully faithful, this implies by 
abstract nonsense that $\mu_X$ and $\eta_X$ are equivalences.
\end{proof}

Denote by $\CatCSt (X)$ the full 2-subcategory of $\CatSt(X)$ of constant 
stacks on $X$. We get

\begin{cor}\label{cstacks}
     If $X$ is 1-connected, the functors 
     $$
     \xymatrix@C=3em{{\CatCat} \ar@<.5ex>[r]^-{(\,\cdot\,)_X} & 
     {\CatCSt(X)} \ar@<.5ex>[l]^-{\catGamma(X, \cdot)} }
     $$
     are equivalences of 2-categories, inverse one to each other.
\end{cor}

\begin{thm}\label{muequiv}
     Let $X$ be locally relatively 1-connected \footnote{Recall that a 
     topological space $X$ is locally relatively 1-connected if each point 
     $x\in X$ has a fundamental system of arcwise connected neighborhoods $U$ 
     such that each loop $\gamma$ in $U$ is homotopic in $X$ to a constant 
     path.}. Then the monodromy
     $$\mu_X \colon \catGamma_X \lra \cat{Y}_{\Pi_1(X)}$$ 
     is an equivalence of 2-functors.
\end{thm}

\begin{proof}
Fix a category $\catc$. By Propositions \ref{1-faith} and \ref{1-full}, 
it remains to show that $\mu_{X,\catc}$ is essentially surjective.
Let us first suppose that $\catc$ is complete, hence we can work in 
the category of sheaves with values in $\catc$.\\
Consider a functor $\alpha\colon \Pi_1(X)\lra\catc$. Define
$$ 
\shv = \Big\{(V,x)\ |\ x\in V,\, V \text{ relatively 
1-connected open subset of }X \Big\} 
$$
and set $(V,x)\leqs (W,y)$ if and only if $W\subset V$, which turns $\shv$
into a category.\\
Let $U\subset X$ be an open subset. We set
$$ \shf_\alpha(U)=\varprojlim_{(V,x)\in\shv \atop V\subset U} \alpha(x) $$
where for any $(V,x)\leqs (W,y)$ we chose a path $\gamma_{xy}\colon x\ra y$ 
in $V$ and use the isomorphism $\alpha(\gamma_{xy})\colon \alpha(x)\ra 
\alpha(y)$ in the projective system. Note that since $V$ is relatively 
1-connected, this automorphism does not depend on the choice of 
$\gamma_{xy}$.\\
Now let $U=\bigcup_{i\in I}U_i$ be a covering stable by finite intersection. 
In order to prove that $\shf_{\alpha}$ is a sheaf, we have to show that the 
natural morphism
$$ 
\shf_{\alpha}(U)=\varprojlim_{(V,x)\in\shv \atop V\subset U} \alpha(x) \to
\varprojlim_{i\in I}\varprojlim_{(V,x)\in\shv \atop V\subset U_i} \alpha(x)=
\varprojlim_{i\in I} \shf_\alpha(U_i)
$$
is an isomorphism. To define the inverse isomorphism, note that 
for $(W,y)\in\shv$ the isomorphisms $\alpha(\gamma_{xy})^{-1}$ define the 
natural isomorphism 
$$ \alpha(y) \isoto \varprojlim_{(i,(V,x))
    \atop V\subset U_i\cap W} \alpha(x). $$
Its inverse defines
$$ \varprojlim_{(i,(V,x)) \atop V\subset U_i} \alpha(x)  \to
 \varprojlim_{(i,(V,x)) \atop V\subset U_i\cap W} \alpha(x)
 \simeq \alpha(y), $$ 
which is compatible with the projective system.\\
By construction it is clear that, if $U$ is relatively 1-connected, then for 
any choice of $x\in U$ we get a natural isomorphism $\shf_\alpha(U)\simeq
\alpha(x)$ (this isomorphism being compatible with restrictions). Hence, 
since relatively 1-connected open subsets form a base of the topology 
of $X$, we get that $\shf_{\alpha}$ is a locally constant sheaf that is 
constant on every relatively 1-connected open subset of $X$. \\
To calculate the monodromy, consider first a path $\gamma\colon x\ra y$
such that there exists a relatively 1-connected open neighborhood of $\gamma$. 
Obviously we get that $\mu(\shf_\alpha)(\gamma)$ is naturally isomorphic to 
$\alpha(\gamma)$. For a general $\gamma$, we decompose $\gamma$ in a finite 
number of paths that can be covered by relatively 1-connected open subsets to 
get the result.\\
Now the general case. Embed $\catc$ into $\widehat{\catc}$ by the 
Yoneda functor
$$ Y\colon\catc \lra \widehat{\catc} .$$ 
Then, given a representation $\alpha$, we can construct $\shf_\alpha$ as a 
locally constant sheaf with values in $\widehat{\catc}$. Then $\shf_{\alpha}$ 
has representable stalks and is therefore in the essential image of the fully 
faithful functor
$$ \catGamma(X,\catc_X) \lra \catGamma(X,\widehat{\catc}_X). $$
Since by construction the monodromy of a locally constant sheaf with values in 
$\catc$ can be calculated by considering it as a locally constant sheaf with 
values in $\widehat{\catc}$, we can conclude.
\end{proof}

Thanks to the 2-Yoneda lemma (as stated for example in 
\cite[cap. 1]{Leinster}), we immediatly recovers the following
\begin{cor}
Let $X$ be locally relatively 1-connected. Then there is an equivalence
of categories 
$$\Pi_1(X)\simeq\catHom(\catGamma_X,\cat{Id}_{\CatCat}),$$
given by $x\mapsto \omega_x$.
\end{cor}

\subsection{Degree 1 non abelian cohomology with constant coefficients}

Let $M$ be a (not necessarily commutative) monoid. Denote by $M[1]$ the 
small category with $1$ as single object and $\on{End}_{M[1]}(1)=M$. Note that
if $G$ is a group then $G[1]$ is a groupoid. Then it is easy to check that we 
get fully faithful functors of categories
$$ [1]\colon\catMon \lra \CatCat \qquad \qquad [1]\colon \catGr\lra\Cat{Gr},$$
where we denote by $\catMon$ the category of small monoids and by $\catGr$ 
that of groups. Also note that if $\catg$ is a connected groupoid (that is
$\pi_0(\catg)\simeq 1$), then for each $P\in \Ob \catg$, the inclusion functor
$\Aut[\catg](P)[1]\lra \catg$ is an equivalence.

Consider the category $\catSet(G)$ of right $G$-sets and $G$-linear maps. 
Then  $G[1]$ is equivalent to the full sub-category of $\catSet(G)$ with $G$ 
as single object. Hence the stack $G[1]_X$ is equivalent to the stack 
$\stkTors(G_X)$ of torsors over the sheaf $G_X$\footnote{Recall that 
$\stkTors(G_X)$ is the stack which associates to each open subset $U\subset X$ 
the category $\catTors(G_U)$ of right $G_U$-sheaves locally free of rank one.
Note that $\catTors(G_X)$ is equivalent to the category of $G$-coverings over 
$X$.}. 

Assume that $X$ is locally relatively 1-connected. 
By Theorem \ref{muequiv} there is an equivalence of categories 
$$\catTors(G_X) \overset{\sim}{\lra} \catHom(\Pi_1(X),G[1]).$$ 
A standard cocycle argument shows that there is an isomorphism of pointed 
sets
$$
\pi_0(\catTors(G_X))\isoto  H^1 (X;G_X).
$$
Assume moreover that the space $X$ is connected. Let us calculate the set 
pointed $\pi_0(\catHom(\Pi_1(X),G[1]))$.
Since $\Pi_1(X)$ is connected, $\Pi_1(X)$ is equivalent to $\pi_1(X)[1]$ for a
choice of a base point in $X$. Hence there is a natural surjective map
$$
\Hom[\catGr](\pi_1(X),G)\ra \pi_0(\catHom(\pi_1(X)[1],G[1])),
$$
and one checks immediately that two morphisms of groups 
$\varphi,\psi\colon \pi_1(X)\to G$ give isomorphic functors if and only if 
there exists $g\in G$ such that $\varphi=\ad (g)\circ\psi$, where $\ad (g)$ is 
the group automorphism of $G$ given by $h\mapsto ghg^{-1}$ for each $h\in G$ 
(automorphisms of this form are called inner automorphisms of $G$).  Hence 
$$
\Hom[\catGr](\pi_1(X),G)/G\simeq \pi_0(\catHom(\pi_1(X)[1],G[1])),
$$
where $G$ acts on the left on $\Hom[\catGr](\pi_1(X),G)$ by conjugation.
We get the classical

\begin{prop}[Hurewicz's formula]\label{H^1nonab}
     Let $X$ be connected and locally relatively 1-connected. 
     Then for any group $G$ there is an isomorphism of pointed 
     sets
     $$  H^1 (X;G_X)\isoto\on{Hom}_{\catGr}(\pi_1(X),G)/G.$$
\end{prop}
\noindent
In particular, if $G$ is commutative one recovers the isomorphism of groups
$$  H^1 (X;G_X)\isoto\on{Hom}_{\catGr}(\pi_1(X),G). $$

More generally, to each complex of groups $G^{-1}\to[d]G^0$ one associates
a small groupoid, which we denote by the same symbol, as follows: 
objects are the elements $g\in G^0$ and morphisms $g\to g'$ are given by 
$h\in G^{-1}$ such that $d(h) g = g'$. If moreover $G^{-1}\to[d]G^0$ has 
the structure of crossed module\footnote{Recall that a complex of groups 
$G^{-1}\to[d]G^0$ is a crossed module if it is endowed with a (left) action of 
$G^0$ on $G^{-1}$ such that (i) $d({}^g h)=\ad (g)(d(h))$ and (ii) 
${}^{d(\tilde h)}h = \ad (\tilde h)(h)$ for any $h,\tilde h \in G^{-1}$ and 
$g\in G^0$. We use the convention as in \cite{Breen} for which $G^i$ is in 
$i$-th degree.}, the associated category is a strict $gr$-category, 
\emph{i.e.} a group object in the category of groupoids. In fact, all strict 
$gr$-categories arise in this way (see for example \cite{Brown-Spencer}, and 
\cite[expos\'e XVIII]{SGA4} for the commutative case). In particular, if $G$ 
is a group, the groupoid $G[1]$ is identified with $G\to 1$ and it has the 
structure of a strict $gr$-category if and only if $G$ is commutative. 
Moreover, the strict $gr$-category $\catAut(G[1])$ of auto-equivalences of 
$G[1]$ is equivalent to $G\to[\ad]\Aut[\catGr](G)$.

Consider the constant stack $(G\to[\ad]\Aut[\catGr](G))_X$. It is 
equivalent to the $gr$-stack $\stkBitors(G_X)$ of $G_X$-bitorsors, 
\emph{i.e.} $G_X$-torsors with an additional compatible structure of left 
$G_X$-torsors (see \cite{Breen} for more details).
Suppose that $X$ is locally relatively 1-connected. 
Then, by Theorem \ref{muequiv}, there is an equivalence of
$gr$-categories 
$$
\catBitors(G_X)\overset{\sim}{\lra} \catHom(\Pi_1(X),G\to[\ad]\Aut[\catGr](G)).
$$ 
One may show (see \emph{loc. cit.}) that 
$$
\pi_0(\catBitors(G_X))\simeq  H^0 (X;G_X\to[\ad]\Aut[\catGr](G_X))
$$
where the right hand side is the 0-th (hyper-)cohomology group of $X$ with 
values in the sheaf of crossed modules $G_X\to[\ad]\Aut[\catGr](G_X)$. 
Suppose that $X$ is connected. Then $\Pi_1(X)\simeq\pi_1(X)[1]$ 
and a similar calculation as above leads to the isomorphism of groups
$$
\pi_0(\catHom(\pi_1(X)[1],G\to[\ad]\Aut[\catGr](G)))\simeq
 \Hom[\catGr](\pi_1(X),\on{Z}(G))\rtimes\on{Out}_{\catGr}(G),
$$
where $\on{Z}(G)$ denotes the center of $G$ and $\on{Out}_{\catGr}(G)=
\Aut[\catGr](G)/G$ is the group of outer automorphisms of $G$, 
which acts on the left on $\Hom[\catGr](\pi_1(X),\on{Z}(G))$ by composition. 
We get

\begin{prop}[Hurewicz's formula II]
     Let $X$ be connected and locally relatively 1-connected. 
     Then for any group $G$ there is an isomorphism of groups
     $$
     H^0 (X;G_X\to[\ad]\Aut[\catGr](G_X))\isoto 
     \Hom[\catGr](\pi_1(X),\on{Z}(G))\rtimes\on{Out}_{\catGr}(G).
     $$
\end{prop}

\noindent
A similar result holds replacing $G\to[\ad]\Aut[\catGr](G)$ by a 
general crossed module $G^{-1}\to[d]G^0$. More precisely, noticing that 
$\ker d$ is central in $G^{-1}$, one gets an isomorphism of groups
$$
H^0 (X;G^{-1}_X\to[d]G^0_X)\isoto \Hom[\catGr](\pi_1(X),\ker d)\rtimes\coker d.
$$

\section{Classification of locally constant stacks}\label{locstacks}

Following our presentation of 1-monodromy, we will approach the theory of 
2-monodromy in the setting of 2-stacks. We refer to \cite{Breen} for the basic 
definitions. Let $X$ be a topological space and let $\TwoCat$, $\TwoPSt(X)$ 
and $\TwoSt(X)$ denote the 3-category{\begin{footnote}{Here and in the 
following, we will use the terminology of 3-categories and 3-functors only 
in the framework of strict 3-categories, \emph{i.e.} categories enriched in 
$\TwoCat$.}\end{footnote}} of small 2-categories, of 2-prestacks and that of 
2-stacks on $X$, respectively. As for sheaves and stacks, there exists a 
2-stack associated to a 2-prestack:

\begin{prop}\label{3-adj}
     The forgetful 3-functor
     $$ \Cat{For}\colon \TwoSt(X) \lra \TwoPSt(X) $$
     has a left adjoint 3-functor
     $$ \ddagger\colon \TwoPSt(X) \lra \TwoSt(X). $$
\end{prop}

\noindent
Hence, we may associate to any 2-category a constant 2-prestack on $X$ and set

\begin{defi}
     Let $\Catc$ be a 2-category. The constant 2-stack on $X$ with stalk 
     $\Catc$ is the image of $\Catc$ by the 3-functor
     $$ 
     (\,\cdot\,)_X\colon \TwoCat \lra \TwoPSt(X) \overset{\ddagger}{\lra} 
     \TwoSt(X).
     $$
     An object $\stks\in\Ob{\Catc_X}(X)$ is called a locally constant stack 
     on $X$ with values in $\Catc$. A locally constant stack is constant 
     with stalk $\on{P}$, if it is equivalent to 
     $\Two{\eta}_{X,\Catc}(\on{P})$ for some object $\on{P}\in\Ob{\Catc}$, 
     where the 3-functor
     $$ \Two{\eta}_{X,\Catc}\colon\Catc \lra \Catc_X(X)  $$
     is induced by the 3-adjunction of Proposition \ref{3-adj}.
\end{defi}

Let us look at the case when $\Catc=\CatCat$, the 2-category of all small
categories. It is easy to see that a stack $\stks$ on $X$ is locally constant 
if and only if there exists an open covering $X=\bigcup {U_i}$ such that 
$\stks|_{U_i}$ is equivalent to a constant stack (as defined in the first 
part). We denote by $\CatLcSt(X)$ the full 2-subcategory of $\CatSt(X)$ whose 
objects are the locally constant stacks.
Similarly, suppose that $\Catc$ is a 2-category that admits all small 2-limits.
Then one can define the notion of a stack with values in $\Catc$ similarly to 
the case of sheaves (see Appendix B). It is not difficult to see that 
the sub-2-category of stacks with values in $\Catc$ which are locally constant 
is naturally equivalent to $\Catc_X(X)$.
For a more general $\Catc$, we can always embed $\Catc$ by the 2-Yoneda lemma 
into the strict 2-category $\widehat \Catc = \CatHom (\Catc^{\op}, \CatCat )$, 
which admits all small 2-limits. Then $\Catc_X(X)$ can be identified with the 
(essentially) full sub-2-category of $\widehat\Catc_X(X)$ defined by objects 
whose stalk is 2-representable.\\

We now follow Section \ref{locsheaves} step by step to define the 2-monodromy 
2-functor. 

Let $f\colon X \ra Y$ be a continous map. We leave to the reader to define
the 3-adjoint 3-functors $f_*$ and $f^{-1}$. 
\begin{defi}
     Denote by $\Cat{\Gamma}(X,\,\cdot\,)$ the 3-functor of global sections
     $$\TwoSt(X)\lra\TwoCat \quad ; \quad \twostks\mapsto 
     \Cat{\Gamma}(X,\twostks)=\twostks(X) $$
     and set $\Cat{\Gamma}_X=\Cat{\Gamma}(X,\cdot)\circ (\cdot)_X$. 
\end{defi}

Since the 3-functor
$$\Cat{\Gamma} (\{\pt\},\cdot)\colon\TwoSt(\{\pt\})\lra\TwoCat $$
is an equivalence of 3-categories, then the 3-functor 
$\Cat{\Gamma}(X,\,\cdot\,)$ is right 3-adjoint to $(\,\cdot\,)_X$.\\
It is not hard to see that we get a 3-transformation of 3-functors $f^{-1}$ 
compatible with $\Two{\eta}_X$ and $\Two{\eta}_Y$, \emph{i.e.} the following 
diagram commutes up to a natural invertible 2-modification:
$$ 
\xymatrix@R=.4cm@C=.5cm{
      {\Cat{\Gamma}_Y} \ar[rr]^{f^{-1}} & & {\Cat{\Gamma}_X} .\\
      & {\Cat{Id}_{\TwoCat} }\ar[ul]^{\Two{\eta}_Y} \ar[ur]_{\Two{\eta}_X} & } 
$$
Similarly to the case of sheaves, this implies that for each point $x\in X$ 
and any locally constant stacks $\stks,\stkt\in \Catc_X(X)$, the natural 
functor
$$ 
i_x^{-1}\stkHom[\Catc_{X}](\stks,\stkt) \lra \catHom[\Catc]
(i^{-1}_x\stks,i_x^{-1}\stkt) 
$$ 
is an equivalence (here $i_x\colon \{\pt\}\to X$ denotes the natural map to 
$x$ and we identify $\Catc$ with global sections of $\Catc_{\pt}$).
Therefore, for each continuous map $f\colon X\ra Y$, we get a natural 
equivalence of locally constant stacks of categories 
$$ 
f^{-1}\stkHom[\Catc_{Y}](\stks,\stkt)  \overset{\sim}{\lra}
 \stkHom[\Catc_X](f^{-1}\stks,f^{-1}\stkt). 
$$

\subsection{The 2-monodromy 2-functor}

Let us prove first the fundamental Lemma of homotopy invariance of 
locally constant stacks.
\begin{lemma}\label{2XxI}
     Consider the maps
     $
     \xymatrix@C=2em{{X} \ar@<.5ex>[r]^-{\iota_t} & 
     {X\times I} \ar@<.5ex>[l]^-{p} }
     $ 
     as in Lemma \ref{XxI}. Then the 3-transformations 
     $$
     \xymatrix@C=3em{ {\Cat{\Gamma}_X} \ar@<.5ex>[r]^-{p^{-1}} 
     & {\Cat{\Gamma}_{X\times I}} \ar@<.5ex>[l]^-{\iota_t^{-1}} }
     $$
     are equivalences of 2-functors, quasi-inverse one to each other. 
\end{lemma}

\begin{proof}
Let $\Catc$ be a 2-category. Since $\iota_t^{-1}\circ p^{-1}\simeq 
\cat{Id}_{\Cat{\Gamma}_X}$, it remains to check that for each 
$\stks\in\Cat{\Gamma}(X\times I,\Catc_{X\times I})$ there is a natural 
equivalence of stacks $p^{-1}\iota_t^{-1}\stks\simeq \stks$.\\
First, let us suppose that $\stks$ is a locally constant stack of categories 
and let us prove that the natural functor
$$\iota_t^{-1}\colon \Gamma(X\times I,\stks)\lra \Gamma(X,\iota_t^{-1}\stks)$$
is an equivalence. Since the sheaves $\shhom[\stks]$ are locally constant, by 
Lemma \ref{XxI} it is clear that this functor is fully faithful. 
Let us show that it is essentially surjective.\\
Note that, since $\stks$ is locally constant, it is easy to see that for 
every open neighborhood $U\times I_j\owns (x,t)$ such that $I_j$ is an 
interval and $\stks|_{U\times I_j}$ is constant, there exists an open subset 
$\tilde U\owns x$ so that for every locally constant sheaf 
$\lshf\in\stks(\tilde U\times I_j)$ there exists $\tilde\lshf\in\stks(\tilde U
\times I)$ such that $\tilde\lshf|_{\tilde U\times I_j}\simeq\lshf$.\\
Now take $\lshf\in\Gamma(X,\iota_t^{-1}\stks)$. Then we can find a covering 
$X\times \{t\}\subset\bigcup_{j\in J}{U_j\times I_j}$ where $I_j$ are open 
intervals containing $t$ such that $\stks|_{U_j\times I_j}$ is constant and 
we can find objects $\lshf_j\in \stks(U_j\times I_j)$ such that 
$\iota_t^{-1}\lshf_j\simeq\lshf|_{U_j}$. Then the isomorphism
$$ 
\shhom[\stks](\lshf_i|_{U_{ij}\times I_{ij}},\lshf_j|_{U_{ij}\times I_{ij}})
\isoto\shhom[\stks](\lshf|_{U_i}|_{U_{ij}},\lshf|_{U_j}|_{U_{ij}}) 
$$ 
implies that we may use the descent data of $\lshf$ to patch the $\lshf_i$ to 
a global object on $X\times I$ that is mapped to $\lshf$ by $\iota^{-1}_t$.\\
The rest of the proof is similar to that of Lemma \ref{XxI}. Consider the stack
of functors $\stkHom[\Catc_{X\times I}](p^{-1}\iota_t^{-1}\stks,\stks)$. Since 
$\stks$ is locally constant, the $\stkHom[\Catc_X]$ stack is locally constant 
and the natural functor 
$$
\iota_t^{-1}\stkHom[\Catc_{X\times I}](p^{-1}\iota_t^{-1}\stks,\stks)\lra 
\stkHom[\Catc_X](\iota_t^{-1}\stks,\iota_t^{-1}\stks)
$$
is an equivalence. We have thus shown that the natural functor 
$$ 
\Gamma(X\times I,\stkHom[\Catc_{X\times I}](p^{-1}\iota_t^{-1}\stks,\stks))
\lra \Gamma(X,\stkHom[\Catc_X](\iota_t^{-1}\stks,\iota_t^{-1}\stks)) 
$$
is an equivalence. We can therefore lift the identity of $\iota_t^{-1}\stks$ 
to get an equivalence 
$$p^{-1}\iota_t^{-1}\stks\overset{\sim}{\lra}\stks. $$
\end{proof}

\begin{cor}
     For each 2-category $\Catc$, the 2-functor
     $$ \Two{\eta}_{I,\Catc}\colon \Catc \lra \CatGamma (I,\Catc_I)  $$
     is an equivalence.
\end{cor}

Hence, for any $X$ and any $t\in I$, we have the equivalence 
$\cat{Id}_{\CatGamma_X} = (p\circ \iota_t)^{-1}\simeq \iota_t^{-1}p^{-1}$. 
Then, for any $s,t\in I$, there exists a canonical (\emph{i.e.} unique up to 
unique invertible modification) equivalence $\iota_t^{-1}\simeq \iota_s^{-1}$. 
With patience, one deduces the following technical Lemma:

\begin{lemma}\label{XxI^2}
\begin{itemize}     
     \item[(i)] For any $s,t,t'\in I$, the topological diagram on the left 
     induces  the diagram of equivalences on the right, which is 
     commutative up to natural invertible modification:
     $$ 
     \xymatrix{ X \ar[r]^{\iota_t} \ar[d]_{\iota_s} & X\times I
     \ar[d]^{\iota_s\times\id_I} \\
     X\times I \ar[r]_{\id_X\times\iota_t} & X\times I^2 } 
     \qquad \qquad
     \xymatrix{ \iota_t^{-1}(\iota_s\times\id_I)^{-1} \ar[r]^{\sim}
     \ar[d]_{\wr} & 
     \iota_{t'}^{-1}(\iota_s\times \id_I)^{-1}  \ar[d]^{\wr}  \\
     \iota_s^{-1}(\id_X\times\iota_t)^{-1} \ar[r]_{\sim} & 
     \iota_{s}^{-1}(\id_X\times\iota_{t'})^{-1} .}
     $$

     \item[(ii)] For any $r,s,t,t'\in I$, the topological diagram
     $$ 
     \xymatrix{ & X \ar[rr]^{\iota_r} \ar[dl]_{\iota_s} \ar[dd]_(.3){\iota_t} 
     & & X\times I \ar[dl]^{\iota_s\times\id_I}\ar[dd]^{\iota_t\times \id_I} \\
     X\times I  \ar[rr]_(.7){\id_X\times\iota_r} \ar[dd]_{\id_X\times\iota_t} 
     & & X\times I^2 \ar[dd]^(.3){\id_X\times\iota_t\times \id_I}  & \\
     & X\times I \ar[rr]_(.3){j_r} \ar[dl]_{\iota_s\times \id_I}  & & 
     X\times I^2 \ar[dl]^{\iota_s\times \id_{I^2}} \\
     X\times I^2 \ar[rr]_{\id_{X\times I}\times\iota_r} & & X\times I^3 & } 
     $$
     induces a (very big) commutative diagram of the corresponding 
     modifications.
     
     \item[(iii)]
     Let $f\colon X\to Y$ be a continous map and $H_f\colon X\times I\ra Y$ 
     the constant homotopy of $f$. Then, for any $t,t'\in I$, the diagram
     $$ 
     \xymatrix{  f^{-1} \ar[r]^{\id_{ f^{-1}}} \ar[d]_{\wr} & f^{-1} 
     \ar[d]^{\wr} \\
     \iota_{t}^{-1}H_f^{-1} \ar[r]_{\sim} & \iota_{t'}^{-1}H_f^{-1} }
     $$
     commutes up to natural invertible modification. 
   \end{itemize}
\end{lemma}

Let $\TwoTop$ denote the 3-category whose objects are topological spaces, 
1-arrows are continuous maps, 2-arrows are homotopies between continuous 
maps and 3-arrows are homotopy classes of homotopies (between homotopies of 
maps). Following the same lines of the Proposition \ref{the2functorgamma} and 
using the two Lemma above, one can then prove

\begin{prop}
     The assignment $(\Catc,X)\mapsto \Cat{\Gamma}(X,\Catc_X)$ defines a 
     3-functor
     $$ \Cat{\Gamma}\colon \TwoCat\times \TwoTop^{\op} \lra \TwoCat,$$
     and the natural 2-functors $\Two{\eta}_{X,\catc}$ define a 
     3-transformation
     $\Two{\eta}\colon\Cat{Q}_1 \lra \Cat{\Gamma}.$
\end{prop}

\begin{defi}
     The homotopy 2-groupoid of $X$ is the 2-groupoid\footnote{Recall that a 
     2-groupoid is a 2-category whose 2-arrows are invertible and 1-arrows are 
     invertible up to a 2-arrow.} 
     $$\Pi_2(X)=\CatHom[\TwoTop](\{\pt\},X). $$
\end{defi}
Roughly speaking, its objects are the points of $X$ and, if $x,y\in X$, 
$\catHom[\Pi_2(X)](x,y)$ is the category $\Pi_1(P_{x,y}X)$, where $P_{x,y}X$ 
is the space of paths starting from $x$ and ending in $y$, endowed with the 
compact-open topology. Compositions laws are defined in the obvious way.
Note that, in particular, $\Two{\pi}_0(\Pi_2(X))=\pi_0(X)$ and for each 
$x\in X$, $\Pic_{\Pi_2(X)}(x)=\pi_0(\Pi_1(\Omega_x X))=\pi_1(X,x)$, where
we denote by $\Omega_x X $ the loop space $P_{x,x}X$ with base point $x$, 
and $Z_{\Pi_2(X)}(x)=\pi_2(X,x)$\footnote{Note 
that the categorical action of $\Pic_{\Pi_2(X)}(x)=\pi_1(X,x)$ on 
$Z_{\Pi_2(X)}(x) = \pi_2(X,x)$ is exactly the classical one of algebraic 
topology.}. We refer to \cite{Hardie-Kamps-Kieboom} for an explicit 
construction of a strictification of $\Pi_2(X)$ when $X$ is Hausdorff.

Let $2\Cat{Gr}$ denote the 3-category of 2-groupoids. Then we have a 3-functor 
$$
\Pi_2\colon\TwoTop\lra 2\Cat{Gr},
$$ 
which defines the Yoneda type 3-functor
$$ 
\Cat{Y}_{\Pi_2}\colon\TwoCat \times \TwoTop^{\op} \lra \TwoCat, \quad 
(\Catc,X) \mapsto \CatHom(\Pi_2(X),\Catc) .
$$

\begin{defi}
     The 2-monodromy 3-transformation
     $$\Two{\mu}^2\colon \Cat{\Gamma} \lra \Cat{Y}_{\Pi_2} $$
     is defined in the following manner. For each topological space $X$, the 
     2-functor 
     $$ 
     \Cat{\Gamma}_{X,\Catc}\colon \CatHom[\TwoTop](\{\pt\},X) \lra
     \CatHom(\Cat{\Gamma}(X,\Catc_X),\Catc) 
     $$
     induces by evaluation a natural 2-functor
     $$  \Cat{\Gamma}(X,\Catc_X) \times\Pi_2(X) \lra \Catc, $$
     hence by adjunction a 2-functor
     $$ \Two{\mu}^2_{X,\Catc}\colon \Cat{\Gamma}(X,\Catc_X) \lra 
     \CatHom(\Pi_2(X),\Catc). $$
\end{defi}

As in the case of 1-monodromy, let us visualize this construction using stalks.
To every 2-stack $\twostks$ we can associate its stalk at $x\in X$, which is 
the 2-category 
$$ \twostks_x=\threecolim{x\in U}{\ \twostks(U)}, $$ 
and a natural 2-functor $\Two{\rho}_{x,\twostks}\colon\twostks (X)\lra
\twostks_x$ (in the case that $\twostks$ is the 2-stack $\twostkSt_X$ of all 
stacks on $X$, we can chose $\Two{\rho}_{x,\twostkSt}=\catF_x$, the ordinary 
stalk 2-functor). Hence we get the natural stalk 3-functor 
$$
\Cat{F}_x\colon\TwoSt (X) \lra\TwoCat \quad ; \quad \twostks\mapsto\twostks_x 
$$
which induces an equivalence 
$$ \Two{\stkLoc}_x \overset{\sim}{\lra} \CatCat $$
(if $\Catc$ is a 2-category, then $(\Catc_X)_x\simeq \Catc$). 
Then similarly as in the case of sheaves, one proves the following

\begin{prop}
     The assignment $(\Catc,(X,x))\mapsto (\Catc_X)_x$ defines a 3-functor
     $$ \Cat{F}\colon \TwoCat\times \TwoTop_*^{\op} \lra \TwoCat $$
     and the 2-functors $\Two{\rho}_{x,\Catc}$ define a 3-transformation
     $ \Two{\rho}\colon \Cat{\Gamma} \lra \Cat{F}. $ 
\end{prop}

We find that
$$ \Two{\rho}\circ \Two{\eta}\colon\Cat{Q_1} \lra \Cat{F} $$
is an equivalence of 3-functors. Let $\Two{\varepsilon}\colon\Cat{F}\lra
\Cat{Q}_1$ be a fixed quasi-inverse to $\Two{\rho}\circ \Two{\eta}$ and set 
$\Two{\omega}=\Two{\varepsilon}\circ \Two{\rho}$.
Fix a topological space $X$ and a locally constant stack $\stks\in 
\CatGamma(X,\Catc_X)$. Then (up to a natural equivalence)
$$ \Two{\mu}^2_{X,\Catc}(\stks)(x) = \Two{\omega}_{x,\Catc}(\stks)$$
(if $\stks$ is a locally constant stack with values in $\CatCat$, then 
$\Two{\omega}_{x,\CatCat}(\stks)$ can be canonically identified with 
$\stks_x$). If $\gamma\colon x\to y$ is a path, then the equivalence 
$\Two{\mu}^2_{X,\Catc}(\stks)(\gamma)$ is defined by the chain of equivalences
$$ 
\Two{\omega}_{x,\Catc}(\stks) \simeq \Two{\omega}_{0,\Catc}(\gamma^{-1}\stks) 
\simeq \Two{\eta}_{I,\Catc}(\gamma^{-1}\stks)  \simeq 
\Two{\omega}_{1,\Catc}(\gamma^{-1}\stks) \simeq\Two{\omega}_{y,\Catc}(\stks), 
$$
where $\Two{\eta}_{I,\Catc}$ is just the global section functor in the case of
ordinary stacks, \emph{i.e.} for $\Catc=\CatCat$. 
If $H\colon\gamma_0\to\gamma_1$ is an homotopy, then the invertible 
transformation $\Two{\mu}^2_{X,\Catc}(\stks)(H)$ is defined by the diagram of 
equivalences
$$ 
\xymatrix{ {\Two{\omega}_{(0,0)}(H^{-1}\stks)}
    \ar[rr]^{\Two{\mu}^2(\stks)(H(\cdot,0))}  
    \ar[dd]_{\Two{\mu}^2(\stks)(\gamma_0)}   & & 
    {\Two{\omega}_{(1,0)}(H^{-1}\stks)}
    \ar[dd]^{\Two{\mu}^2(\stks)(\gamma_1)}\\
    &{\Two{\eta}_{I^2}(H^{-1}\stks)}\ar[dr]_{\rho_{(1,1)}}
    \ar[dl]^{\rho_{(0,1)}} 
    \ar[ul]_{\rho_{(0,0)}} \ar[ur]^{\rho_{(1,0)}} & \\
    {\Two{\omega}_{(0,1)}(H^{-1}\stks)}
    \ar[rr]_{\Two{\mu}^2(\stks)(H(\cdot,1))}
    && {\Two{\omega}_{(1,1)}(H^{-1}\stks)}.} 
$$
In particular, the following diagram of 
3-transformations commutes (up to a natural 2-modification)
$$ 
     \xymatrix@C=.7cm@R=.5cm{ {\Cat{\Gamma}}\ar[rr]^-{\Two{\mu}^2}  
     \ar[dr]_-{\Two{\omega}} && {\Cat{Y}_{\Pi_2},} \ar[dl]^-{\cat{ev}} \\
     &  {\Cat{Q}_1} & }
$$
where $\cat{ev}$ denotes the evaluation 3-transformation.

\subsection{Classifying locally constant stacks}

Let $\cat{\Delta}\colon\Cat{Q}_1\lra\Cat{Y}_{\Pi_2}$ denote the diagonal 
3-transformation. Exactly as in the sheaf case, one gets
\begin{prop} 
     The image of a constant stack is equivalent to a trivial representation.
     In other words, the diagram of 3-transformations 
     $$ 
     \xymatrix@C=.7cm@R=.5cm{ {\Cat{\Gamma}}\ar[rr]^-{\Two{\mu}^2} & & 
     {\Cat{Y}_{\Pi_2}}\\
     & \ar[ul]^-{\Two{\eta}} {\Cat{Q}_1} \ar[ur]_-{\cat{\Delta}} 
     & }     $$
     commutes up to 2-modifications.
\end{prop}

\begin{prop}
     For any topological space $X$ and 2-category $\Catc$, the 2-functor 
     $$
     \Two{\mu}^2_{X,\Catc}\colon \Cat{\Gamma}(X,\Catc_X) 
     \lra\CatHom(\Pi_2(X),\Catc)
     $$ 
     is faithful and conservative.
\end{prop}

\begin{proof}
We have to show that, for any locally constant stack $\stks$ and $\stks'$, the 
induced functor
$$ 
\Two{\mu}^2\colon \catHom[\Catc_X](\stks,\stks') \lra 
\catHom[\CatHom(\Pi_2(X),\Catc)](\Two{\mu}^2(\stks),
\Two{\mu}^2(\stks')) 
$$ 
is faithful and conservative.
Let $F,G\colon\stks\lra \stks'$ be two functors of stacks. 
Since for each $x\in X$, there is a natural isomorphism 
$\Hom(\Two{\mu}^2 (F)(x),\Two{\mu}^2 (F)(x))\simeq \Hom(F_x,G_x),$  
we get the commutative diagram
$$ 
\xymatrix { {\Hom(F,G)} \ar[d]_-{\catF_x} \ar[r]^-{\Two{\mu}^2} & 
          {\Hom(\Two{\mu}^2 (F),\Two{\mu}^2 (G))} \ar[d]^-{ev_x} \\
          {\shhom(F,G)_x} \ar[r]^-{\sim} &  
          {\Hom(F_x,G_x)} .} 
$$
Let $\varphi,\psi\colon F\to G$ be two morphisms of functors. 
If $\Two{\mu}^2 (\varphi) = \Two{\mu}^2 (\psi)$, then $\varphi_x  = \psi_x$ 
for all $x\in X$ and, since $\shhom (F,G)$ is a sheaf, we get $\varphi = \psi$.
Similarly, if $\Two{\mu}^2 (\varphi)$ is an isomorphism, it follows that the 
morphism $\varphi$ is an isomorphism. 
\end{proof}
 
\begin{prop}
     Let $X$ be locally 1-connected. Then for each 2-category $\Catc$, 
     the 2-functor 
     $$
     \Two{\mu}^2_{X,\Catc}\colon \Cat{\Gamma}(X,\Catc_X)
     \lra\CatHom(\Pi_2(X),\Catc)
     $$ 
     is full.
\end{prop}

\begin{proof}
We have to show that the induced functor
$$
\Two{\mu}^2\colon \catHom[\Catc_X](\stks,\stks') \lra 
\catHom[\CatHom(\Pi_2(X),\Catc)](\Two{\mu}^2(\stks),
\Two{\mu}^2(\stks')) 
$$ 
is full and essentially surjective. \\
A morphism $\phi:\,\Two{\mu}^2(\stks)\ra \Two{\mu}^2(\stks')$ is defined by a 
family of functors 
$$\phi_x:\,\on{ev}_x\Two{\mu}^2(\stks) \lra \on{ev}_x\Two{\mu}^2(\stks') $$
and for any path $\gamma:\,x\ra x'$ a canonical isomorphism
$$ 
\xymatrix{
\on{ev}_x\Two{\mu}^2(\stks) \ar[r]^{\phi_x} \ar[d]_{\Two{\mu}^2(\stks)(\gamma)}
 & \on{ev}_x\Two{\mu}^2(\stks') \ar[d]^{\Two{\mu}^2(\stks')(\gamma)} \\
 \on{ev}_y\Two{\mu}^2(\stks) \ar[r]_{\phi_y} \df[ur]^{\sim} 
 & \on{ev}_y\Two{\mu}^2(\stks')  
} 
$$
Since $\stks,\stks'$ are locally constant, we have an equivalence of categories
$\stkHom[\Catc_X](\stks,\stks')_x\simeq\catHom(\stks_x,\stks'_x)$ and we may 
lift $\phi_x$ to a 1-connected open neighborhoods $U_x$ of $x$, say to a 
functor
$$ \varphi^x:\ \stks|_{U_x} \lra \stks'|_{U_x}. $$
Consider $z\in U_x\cap U_y$ and chose paths $\gamma_{xz}$ from $x$ to $z$ and 
$\gamma_{yz}$ from $y$ to $z$. We get the diagram
$$ 
\xymatrix{
  \stks_z \ar[d]_{\on{ev}_z(\varphi^x)} & \stks_x \ar[d]^{\on{ev}_x(\varphi^x)} \ar[r] \ar[l] & 
  \stks_z \ar[d]^{\on{ev}_z(\varphi^z)} & \stks_y \ar[d]^{\on{ev}_y(\varphi^y)} \ar[r] \ar[l] & 
  \stks_z \ar[d]^{\on{ev}_z(\varphi^y)} \\
   \stks'_z & \stks'_x \ar[r] \ar[l] & \stks'_z & \stks'_y \ar[r] \ar[l] & \stks'_z}
$$
that commutes up to natural isomorphism. Moreover the horizontal lines are 
canonically isomorphic to Identities, hence we get a natural isomorphism 
$\on{ev}_z(\varphi^x)\simeq \on{ev}_z(\varphi^y)$ that we can lift to a small 
neighborhood of $z$. Since $U_x$ and $U_y$ are 1-connected, this lift does not 
depend on the choice of the paths $\gamma_{xz}$ and $\gamma_{yz}$ and the 
isomorphism $\on{ev}_z(\varphi^x)\simeq \on{ev}_z(\varphi^y)$ is canonical in a
neighborhood of $z$, i.e. it satisfies the cocycle condition and can be patched
to an isomorphism
$$ \varphi^x|_{U_{xy}} \simeq \varphi^y|_{U_{xy}}. $$
Clearly these isomorphisms satisfy the cocycle condition, hence we get an 
isomorphism $\varphi:\,\stks \lra \stks'$. By construction 
$\Two{\mu}^2(\varphi)\simeq\phi$.\\
Next consider two functors $\varphi,\psi:\,\stks\lra\stks'$ and a morphism
$f:\,\Two{\mu}^2(\varphi) \ra \Two{\mu}^2(\psi)$. Such a morphism is defined by
a family
$$ f_x:\, \on{ev}_x\Two{\mu}^2(\varphi) \ra  \on{ev}_x\Two{\mu}^2(\psi) $$
with compatibility conditions, hence by a family
$$ f_x:\, \varphi_x \ra \psi_x $$ 
such that for every path $\gamma$ from $x$ to $y$
$$ \xymatrix{
  \Two{\mu}^2\stks'(\gamma)\varphi_x \ar[r]^{f_x} \ar[d]_{\sim} & \Two{\mu}^2\stks'(\gamma)\psi_x \ar[d]^{\sim} \\
  \varphi_y\Two{\mu}^2\stks(\gamma)  \ar[r]_{f_y}& \psi_y\Two{\mu}^2\stks(\gamma) } $$
commutes. Lifting $f_x$ to any arcwise connected neighborhood $U_x$ of $x$, this diagram implies
that $\Two{\mu}^2(f_x)=f_z$ for all $z\in U_z$ and we can show that we can patch the $f_x$
similarly to the case of 1-monodromy. 

\end{proof}

\begin{cor}\label{2-con}
     Let $X$ be 2-connected. Then the 2-functor 
     $$\Two{\eta}_{X,\Catc} \colon \Catc \lra \Catc_X (X)$$
     is a natural equivalence of 2-categories. For $\Catc=\CatCat$, 
     the 2-functor
     $\catGamma(X, \cdot)$ provides a quasi-2-inverse to 
     $\Two{\eta}_{X,\Catc}= (\cdot)_X$.
\end{cor}

\begin{proof}
The proof follows the same lines as the proof of Proposition \ref{1-con}, 
using the fact that $\Pi_2(X)\simeq 1$. If $\Catc=\CatCat$, one may chose 
$\Two{\eta}_{X,\Catc}= (\cdot)_X$, hence $\catF_{x_0}$ gives a quasi-2-inverse 
for any choice of $x_0\in X$. Thanks to the natural 2-transformation 
$\catGamma(X, \cdot)\lra \catF_{x_0}$, the 2-functor $\catGamma(X, \cdot)$ 
provides another quasi-2-inverse.
\end{proof}

\begin{thm}\label{mu2equiv}
     Let $X$ be locally relatively 2-connected \footnote{Recall that a 
     topological space $X$ is locally relatively 2-connected if each point 
     $x\in X$ has a fundamental system of 1-connected neighborhoods $U$ such 
     that every homotopy of a path in $U$ is homotopic to the constant 
     homotopy in $X$.}. 
     Then
     $$ \Two{\mu}^2_X\colon  \Cat{\Gamma}_X \lra \Cat{Y}_{\Pi_2(X)} $$
     is an equivalence of 3-functors.
\end{thm}

\begin{proof}
We have to show that for each 2-category $\Catc$, the 2-functor 
$\Two{\mu}^2_{X,\Catc}$ is essentially surjective.\\
Suppose first that $\Catc$ is 2-complete and let 
$\Two{\alpha}\in \CatHom(\Pi_2(X),\CatCat)$. Set
$$ 
\shv = \Big\{(V,x)\ |\ x\in V,\, V \text{ relatively 2-connected open subset 
of } X \Big\} 
$$
and define  $(V,x)\leqs (W,y)$ if and only if $W\subset V$.\\
Let $U\subset X$ be an open subset. We set
$$ 
\stks_{\Two{\alpha}}(U)=\twolim{(V,x)\in\shv \atop V\subset U}{\Two{\alpha}(x)}
$$
where for any $(V,x)\leqs (W,y)$, we chose a path $\gamma_{xy} \colon x\ra y$ 
in $V$ and use the equivalence 
$\Two{\alpha}(\gamma_{xy})\colon\Two{\alpha}(x)\isoto \Two{\alpha}(y)$ in the 
projective system, and for any $(V,x)\leqs (W,y)\leqs (Z,z)$, we chose a 
homotopy $H_{\gamma_{xy}, \gamma_{yz} ,\gamma_{xz}}\colon \gamma_{xy} 
\gamma_{yz} \to \gamma_{xz}$ in $V$ and use the invertible transformation 
of functors $\Two{\alpha}(H_{\gamma_{xy}, \gamma_{yz} ,\gamma_{xz}})$.
Note that since $V$ is relatively 2-connected, the equivalence
$\Two{\alpha}(\gamma_{xy})$ is unique up to invertible transformation and 
the invertible transformation $\Two{\alpha} (H_{\gamma_{xy}, \gamma_{yz} ,
\gamma_{xz}})$ does not depend on the choice of the homotopy
$H_{\gamma_{xy}, \gamma_{yz} ,\gamma_{xz}}$.\\
One argues as in the proof of Theorem~\ref{muequiv} to show that the pre-stack 
defined by $X\supset U\mapsto  \stks_{\Two{\alpha}}(U)$ is actually a stack.
By definition it is clear that, if $U$ is relatively 2-connected, then for any
choice of $x\in U$ (and paths from $x$ to $y$ for every $y\in U$) we get an 
equivalence of categories $\stks_{\Two{\alpha}}(U)\simeq \Two{\alpha}(x)$ 
compatible with restriction functors in a natural sense. 
Hence, the stack $\stks_{\Two{\alpha}}$ is constant on every relatively 
2-connected open subset of $X$. Since relatively 2-connected open subsets form 
a base of the topology of $X$, we get that $\stks_{\Two{\alpha}}$ is locally 
constant.\\ 
The computation of the 2-monodromy of $\stks_{\Two{\alpha}}$ is similar to 
that of 1-monodromy in the proof of Theorem~\ref{muequiv}.\\
For a general 2-category $\Catc$, we can use the 2-Yoneda lemma to reduce to
this last case. 
\end{proof}

Suppose that $X$ is connected and locally relatively 2-connected, and let 
$\Omega X$ be the loop space at a fixed base point $x_0\in X$.
Consider the following diagram of topological space and continuous maps
$$
\xymatrix@C=3.7em{ {(\Omega X)^3 } \ar@<+1.5ex>[r]|-{q_{23}}
 \ar@<+0.5ex>[r]|-{m\times id} \ar@<-0.5ex>[r]|-{id\times m}
 \ar@<-1.5ex>[r]|-{q_{12}} & {(\Omega X)^2 } \ar@<-1ex>[r]|-{q_1}
 \ar[r]|-{m} \ar@<1ex>[r]|-{q_2} & {\Omega X } \ar[r] & \{x_0\} }
$$
where the $q_i$'s, the $q_{ij}$'s and the $q_{ijk}$'s are the natural 
projections and $m$ the composition of paths in $\Omega X$\footnote{Note that 
$\Omega X$ does not define a simplicial topological space, since the maps 
$m\circ (id \times m)$ and $m\circ (m \times id)$ are not equal but only 
homotopic. What one gets is a 2-simplicial object in the 2-category $\CatTop$. 
This will not cause particular difficulties, since for locally constant 
objects there is a natural invertible transformation of functors 
$(m \times id)^{-1}m^{-1} \isoto (id \times m)^{-1}m^{-1}$.}.

Let $\stks$ be a locally constant stack on $X$ with values in $\Catc$.
Theorem \ref{mu2equiv} asserts that $\stks$ is completely and uniquely 
(up to equivalence) determined by its 2-monodromy $\Two{\mu}^2_{X,\Catc}(\stks)
\colon \Pi_2(X)\lra \Catc$. 
Since $X$ is connected, the stalks of $\stks$ are all equivalent. Let us 
denote by $\on{P}$ the stalk at $x_0$. By chosing paths from $x_0$ to any point
$x$, the 2-monodromy reads as a monoidal functor 
$\Two{\mu}^2_{X,\Catc}(\stks)\colon\Pi_1(\Omega X )
\lra \catAut[\Catc](\on{P})$. Since the topological space $\Omega X $ 
satisfies the hypothesis of Theorem \ref{muequiv}, there is a chain 
of equivalences of categories
$$
\catHom(\Pi_1(\Omega X ),\catAut[\Catc](\on{P}))\underset{\mu}{\overset{\sim}
{\longleftarrow}}
\catGamma(\Omega X ,\catAut[\Catc](\on{P})_{\Omega X })\simeq
\catAut[\Catc_X](\Two{\eta}_{\Omega X,\Catc}(\on{P})).
$$
Then the 2-monodromy is equivalent to a pair $(\alpha ,\nu)$ where 
$\alpha\colon\Two{\eta}_{\Omega X ,\Catc}(\on{P}) \overset{\sim}{\lra}
\Two{\eta}_{\Omega X ,\Catc}(\on{P})$ is an equivalence of constant stacks on 
$\Omega X $ and
$$
\nu \colon q_1^{-1}\alpha \circ q_2^{-1}\alpha\overset{\sim}{\ra} m^{-1}\alpha
$$
is an invertible transformation of functors of stacks on $(\Omega X)^2$ such 
that the following diagram of invertible transformations of functors of stacks 
on $(\Omega X)^3$ commutes
\begin{equation}\label{nudiagram}
\xymatrix@C=-2em{ && q_{1}^{-1}\alpha\circ q_{2}^{-1}\alpha
\circ q_{3}^{-1}\alpha \ar[dll]_{\sim} \ar[drr]^{\sim} && \\
q_{12}^{-1}(q_{1}^{-1}\alpha\circ q_{2}^{-1}\alpha)\circ q_{3}^{-1}\alpha
\ar[d]_{\nu} &&&& q_{1}^{-1}\alpha\circ q_{23}^{-1}(q_{1}^{-1}\alpha\circ
q_{2}^{-1}\alpha) \ar[d]^{\nu} \\
q_{12}^{-1}m^{-1}\alpha\circ q_{3}^{-1}\alpha \ar[d]_{\wr} &&&& 
q_{1}^{-1}\alpha\circ q_{23}^{-1}m^{-1}\alpha \ar[d]^{\wr} \\
(m \times id)^{-1}(q_{1}^{-1}\alpha\circ q_{2}^{-1}\alpha) \ar[dr]_{\nu} &&&&
(id\times m)^{-1}(q_{1}^{-1}\alpha\circ q_{2}^{-1}\alpha) \ar[dl]^{\nu} \\
& (m \times id)^{-1}m^{-1}\alpha \ar[rr]_{\sim} && 
(id\times m)^{-1}m^{-1}\alpha \, .&}
\end{equation}
Roughly speaking, $\nu$ is given by a family of functorial invertible 
transformations $\nu_{12}\colon\alpha_{\gamma_1}\circ\alpha_{\gamma_2}\isoto
\alpha_{\gamma_1\gamma_2}$ for any $\gamma_1, \gamma_2\in \Omega X$, 
such that for $\gamma_1,\gamma_2,\gamma_3\in \Omega X$ the following diagram 
commutes
$$
\xymatrix@C=4em{ \alpha_{\gamma_1 \gamma_2}\circ\alpha_{\gamma_3} 
\ar[d]_{\nu_{12,3}} & 
\alpha_{\gamma_1}\circ\alpha_{\gamma_2}\circ\alpha_{\gamma_3} 
\ar[l]_{\nu_{12}} \ar[r]^{\nu_{23}} & 
\alpha_{\gamma_1}\circ\alpha_{\gamma_2 \gamma_3} \ar[d]^{\nu_{1,23}} \\
\alpha_{(\gamma_1 \gamma_2) \gamma_3} \ar[rr]^{\sim} && 
\alpha_{\gamma_1 (\gamma_2 \gamma_3)} .}
$$

\begin{defi}
We call the triplet $(\on{P},\alpha,\nu)$ a descent datum for the locally 
constant stack $\stks$ on $X$.
\end{defi}

Let us analize a particular case, for which the descent datum admits a more 
familiar description. Let $G$ be a (not necesserly commutative) group.
Recall that a $G_X$-gerbe is a stack locally equivalent to the stack of 
torsors $\stkTors(G_X)\simeq G[1]_X$, that is a locally constant stack 
on $X$ with stalk the groupoid $G[1]$. 
Let $\stkg$ be a $G_X$-gerbe and assume 
that $X$ is connected and locally relatively 2-connected. By Morita theorem 
for torsors (see \cite{Giraud} cap. IV), an equivalence 
$\alpha\colon \stkTors (G_{\Omega X})\overset{\sim}{\lra}
\stkTors (G_{\Omega X })$ is given by $\shn\mapsto \shp\wedge\shn$ for a 
$G_{\Omega X}$-bitorsor $\shp$, where $\cdot\wedge\cdot$ denotes the 
contracted product. Then, the descent datum for $\stkg$ is equivalent to the 
datum of $(G[1], \shp, \nu)$ where 
$$\nu \colon q_1^{-1}\shp \wedge q_2^{-1}\shp\overset{\sim}{\ra} m^{-1}\shp$$ 
is an isomorphism of $G_{\Omega X}$-bitorsors on $(\Omega X)^2$ satisfying a 
commutative contraint similar to that of diagram \eqref{nudiagram}.
See \cite{Brylinski} cap. 6, for related constructions of line bundles on the 
free loop space of a manifold.

\subsection{Degree 2 non abelian cohomology with constant coefficients}

Let $\catd$ be a monoidal category. Denote by $\catd [1]$ the 2-category 
with $1$ as single object and $\catEnd[{\catd[1]}](1)=\catd$. Note that if 
$\catd$ is a groupoid whose monoidal structure is rigid\footnote{Recall that a 
monoidal category $(\catd,\otimes)$ is rigid if for each $P\in \Ob\catd$ there 
exists an object $P^*$ and natural ismorphisms $P\otimes P^*\simeq I$ and 
$P^*\otimes P\simeq I$, where $I$ denotes the unit object in $\catd$.}, then 
$\catd [1]$ is a 2-groupoid. It is easy to see that we get a fully faithful 
2-functor 
$$ [1]\colon\CatMon \lra \TwoCat, $$
where we denote by $\CatMon$ the strict 2-category of small monoidal 
categories with monoidal functors and monoidal transformations. This functor 
sends rigid monoidal groupoids to 2-groupoids. We follow the terminology of 
\cite{Breen2} and call $gr$-category a rigid monoidal groupoid. 

Note that if $M$ is a monoid, then $M[1]$ is monoidal if and only if $M$ is 
commutative, and that if $M$ is also a group, then $M[1]$ is a $gr$-category. 
Hence we get fully faithful functors of catgories
$$ 
[2]=[1]\circ [1]\colon \catMonCom \lra \CatMon \lra \TwoCat \qquad , \qquad   
[2]\colon \catGrCom \lra \TwoGr 
$$
where the uppercase $\cat{c}$ means commutative structures. 
Conversely, if $\Cat{G}$ is a connected 2-groupoid, for each object 
$\on{P}\in \Ob\Cat{G}$, the inclusion 2-functor $\catAut[\Cat{G}](\on{P})[1] 
\lra \Cat{G}$ 
is a 2-equivalence. If $\Cat{G}$ is even 1-connected 
(\emph{i.e.} moreover $\catAut[\Cat{G}](\on{P})$ is a connected 
groupoid for some, hence all, $\on{P}$), then $\on{Z}_{\Cat{G}}(\on{P})[2]
\simeq\Cat{G}$. 

For a not necessarily commutative group $G$, we can consider the
strict $gr$-category $\catAut(G[1])$ which gives rise to the 2-groupoid
$$ G[\![2]\!]=\catAut(G[1])[1] $$
Recall that $\catAut(G[1])$ is equivalent to $G\to[\ad] \Aut[\catGr](G)$. 
Hence if $G$ is commutative, then $\catAut(G[1])$ is completely disconnected 
but only $\id\in \Ob \catAut(G[1])$ is $G$-linear, so we get a monoidal functor
$$ \catAut[G](G[1])[1] \simeq G[2] \lra G[\![2]\!]=\catAut(G[1])[1] $$
that identifies $G[2]$ to a sub-2-category of $G[\![2]\!]$ which has only the 
identity 1-arrow but the same 2-arrows.

Consider an object $\catc$ of the 2-category $\CatCat$ (resp. $\CatCat_G$). 
Then $\catEnd(\catc)[1]$ (resp. $\catEnd[G](\catc)[1]$) is just the full 
sub-2-category of $\CatCat$ (resp. $\CatCat_G$) with the single object 
$\catc$. Hence,  $\catEnd(\catc)[1]_X$ (resp. $\catEnd[G](\catc)[1]_X$) is the 
2-stack of locally constant stacks (resp. $G_X$-linear stacks) on $X$ with 
stalk $\catc$.\\
If $X$ is a  locally relatively 2-connected space, by Theorem \ref{mu2equiv} 
equivalence classes of such stacks are classified by the set
\begin{equation}\label{mu2endC}
\Two{\pi}_0(\CatHom(\Pi_2(X),\catAut(\catc)[1])).
\end{equation}
Assume moreover that $X$ is connected. Then the 2-groupoid $\Pi_2(X)$ is 
connected, hence it is equivalent to $\Pi_1(\Omega X)[1]$ for some base point
$x_0\in X$. Hence there is a natural surjective map
$$
\pi_0(\catHom[\otimes](\Pi_1(\Omega X),\catAut(\catc)))\lra 
\Two{\pi}_0(\CatHom(\Pi_1(\Omega X)[1],\catAut(\catc)[1])),
$$
where $\catHom[\otimes](\cdot,\cdot)$ denotes the category of monoidal 
functors. One checks that, given two monoidal functors $\Phi, \Psi \colon
\Pi_1(\Omega X)\lra \catAut(\catc)$, they give equivalent 2-functors if and 
only if there exists an equivalence $\varphi\colon \catc \overset{\sim}{\lra} 
\catc$ and an invertible monoidal transformation $\alpha\colon \Phi (\cdot) 
\circ \varphi \isoto \varphi\circ \Psi (\cdot)$.
We thus get that the set \eqref{mu2endC} is isomorphic to
\begin{equation}\label{mu2endCbis}
\pi_0(\catHom[\otimes](\Pi_1(\Omega X),\catAut(\catc)))/ \Pic(\catc) ,
\end{equation}
where the group $\Pic(\catc) $ acts by conjugation.
A similar reslut holds, replacing $\CatCat$ by $\CatCat_G$.

The previous classification becomes very simple in the following case:
\begin{prop}
     Let $X$ be connected and locally relatively 2-connected. If $\catc$ is a 
     category with trivial center, then the set of equivalence classes of 
     locally constant stacks on $X$ with stalk $\catc$ is isomorphic to
     $H^1(X;\Pic(\catc)_X)$. 
\end{prop}

\begin{proof}
Suppose for simplicity that $\catc$ is a groupoid. Since $\on{Z}(\catc)\simeq 
1$, the monoidal functor $\catAut (\catc) \overset{\pi_0} {\lra} 
\Pic(\catc)[0]$ is an equivalence, where the group $\Pic(\catc)$ is view as 
a discrete category. Hence it easy to check that
$$
\pi_0(\catHom[\otimes](\Pi_1(\Omega X),\Pic(\catc)[0]))\simeq 
\Hom[\catGr](\pi_1(X), \Pic(\catc)).
$$
Hence $\Two{\pi}_0(\CatHom(\Pi_2(X),\catAut(\catc)[1]))$ is isomorphic to 
$\Hom[\catGr](\pi_1(X), \Pic(\catc))/ \Pic(\catc)$, where $\Pic(\catc)$ acts by
conjugation. It remains to apply Proposition \ref{H^1nonab}.
\end{proof}

Let us analyze more in detail the case of gerbes. We start with the abelian 
case. Let $G$ be a commutative group and take $\catc=G[1]$. Since there is an 
obvious equivalence of strict $gr$-categories $G[1]\simeq \catEnd[G](G[1])$, 
we get that $G[2]$ is just the full sub-2-category of $\CatCat_G$ with the 
single object $G[1]$. Hence the 2-category $\Cat{\Gamma}(X,G[2]_X)$ is 
equivalent to the strict 2-category of abelian\footnote{Recall that an abelian 
$G_X$-gerbe is a $G_X$-linear stack locally $G_X$-equivalent to the 
$G_X$-linear stack of torsors $\stkTors(G_X)\simeq G[1]_X$ (we refer to
\cite{Giraud,Breen}, and \cite{Brylinski} for a complete introduction to 
abelian gerbes).} $G_X$-gerbes $\Cat{Ger}_{ab}(G_X)$.\\
By some cocycle arguments (see for example \cite{Brylinski} cap. IV), one 
shows that there is an isomorphism of groups 
$$\Two{\pi}_0(\Cat{Ger}_{ab}(G_X))\simeq  H^2 (X;G_X).$$ 
Assume that $X$ is locally relatively 2-connected. By Theorem 
\ref{mu2equiv}, there is an equivalence of monoidal 2-categories 
$$\Cat{Ger}_{ab}(G_X) \overset{\sim}{\lra} \CatHom(\Pi_2(X),G[2]).$$ 
Since $\Pic_G(G[1])\simeq 1$, if $X$ is connected \eqref{mu2endCbis} 
gives the group $\pi_0(\catHom[\otimes](\Pi_1(\Omega X),G[1]))$.
Hence we get
\begin{prop}[Hurewicz-Hopf's formula]\label{H^2ab}
     Let $X$ be connected and locally relatively 2-connected. 
     Then for any commutative group $G$ there is an isomorphism of groups
     $$ H^2 (X;G_X)\isoto\pi_0(\catHom[\otimes](\Pi_1(\Omega X),G[1])).$$
\end{prop}

To give an explicit description of the right hand side, let us start by
considering a $gr$-category $\cat{H}$. Recall that there exists an "essentially
exact"\footnote{This means that the monoidal functor $i$ (resp. $\pi_0$) is 
fully faithfull (resp. essentially surjective) and that the essential image of 
$i$ is equivalent to the kernel of $\pi_0$ as monoidal categories.} sequence 
of $gr$-categories
$$
\xymatrix{ 1 \ar[r] & {\on{A}_{\cat{H}}[1]} \ar[r]^-{i} & 
{\cat{H}} \ar[r]^-{\pi_0} &  {\pi_0(\cat{H})[0]} \ar[r] & 1 ,}
$$ 
where $\on{A}_{\cat{H}}$ denotes the commutative group $\Aut[\cat{H}](I)$ of 
automorphims of the unit object, and the group $\pi_0(\cat{H})$ is view as 
discrete category, which acts on $\on{A}_{\cat{H}}$ by conjugation. 
Then, if $\catg$ is another $gr$-category, we get an exact sequence of pointed 
sets
\begin{equation}\label{exgr}
\xymatrix{ 1 \ar[r] & {\pi_0(\catHom[\otimes](\pi_0(\cat{H})[0],\catg))} 
 \ar[r] & {\pi_0(\catHom[\otimes](\cat{H},\catg))} \ar[r] &  
 {\pi_0(\catHom[\otimes](\on{A}_{\cat{H}}[1],\catg))} .}
\end{equation}

\begin{lemma}\label{abcohom}
     Let $\catg=G[1]$, for $G$ an abelian group. Then \eqref{exgr} gives an 
     exact sequence of abelian groups
     $$
     \xymatrix{ 1 \ar[r] & {H^2(\pi_0(\cat{H});G)} \ar[r] & 
     {\pi_0(\catHom[\otimes](\cat{H},G[1]))} \ar[r] &  
     {\Hom[\catGr](\on{A}_{\cat{H}},G)} ,}
     $$
     where $G$ is view as a $\pi_0(\cat{H})$-module with trivial action.
\end{lemma} 

\begin{proof}
Set $H=\pi_0(\cat{H})$. It is easy to see that a monoidal 
functor $H[0]\lra G[1]$ is given by a set-theoretic function 
$\lambda \colon H\times H\ra G$ such that
$$ 
\lambda (h_1,h_2)\lambda(h_1h_2,h_3) = \lambda (h_2,h_3)\lambda (h_1,h_2h_3), 
$$
and that two monoidal functors $\lambda,\lambda'$ are isomorphic if and only 
if there exists a function $\nu \colon H\ra G$ such that
$$\lambda(h_1,h_2)\nu(h_1h_2) = \lambda'(h_1,h_2)\nu(h_1)\nu(h_2).$$
Hence we get an isomorphism of groups
$$ \pi_0(\catHom[\otimes](H[0],G[1])\simeq H^2(H;G), $$
where $G$ is view as a $H$-module with trivial action. \\
Similarly, one easily checks that $\pi_0(\catHom[\otimes]
(\on{A}_{\cat{H}}[1],G[1]))$ is isomorphic to 
$\Hom[\catGr](\on{A}_{\cat{H}},G)$. 
\end{proof}

Let $\Hom[\pi_0(\cat{H})](\on{A}_{\cat{H}},G)$ denote the subgroup of 
morphisms in $\Hom[\catGr](\on{A}_{\cat{H}},G)$ which are 
$\pi_0(\cat{H})$-linear. One easily checks that the morphism 
$\pi_0(\catHom[\otimes](\cat{H},G[1])) \to \Hom[\catGr](\on{A}_{\cat{H}},G)$
factors through $\Hom[\pi_0(\cat{H})](\on{A}_{\cat{H}},G)$. Then one gets an 
exact sequence of abelian groups
$$
     \xymatrix{ 1 \ar[r] & {H^2(\pi_0(\cat{H});G)} \ar[r] & 
     {\pi_0(\catHom[\otimes](\cat{H},G[1]))} \ar[r] &  
     {\Hom[\pi_0(\cat{H})](\on{A}_{\cat{H}},G)} \ar[r]^-{\delta} & 
     {H^3(\pi_0(\cat{H});G)} }
$$
where the coboundary morphism $\delta$ is described as follows.\\
Recall that to the $gr$-category $\cat{H}$ one associates a cohomology 
class\footnote{If $\cat{H}= H^{-1}\to[d] H^0$, this class coincides with the 
usual one in $H^3(\coker d; \ker d)$ attached to the crossed module 
$H^{-1}\to[d] H^0$. See for example \cite[cap. V]{Brown}.} $[\cat{H}]$ in 
$H^3(\pi_0(\cat{H});\on{A}_{\cat{H}})$, where $\on{A}_{\cat{H}}$ is endowed 
with the conjugation action of $\pi_0(\cat{H})$ (see for example 
\cite{Breen2}).
Hence, to each $\pi_0(\cat{H})$-linear morphism $f\colon\on{A}_{\cat{H}}\to G$,
one associates the image of $[\cat{H}]$ by the induced morphism 
$\hat f \colon H^3(\pi_0(\cat{H});\on{A}_{\cat{H}})\to H^3(\pi_0(\cat{H});G)$.

\begin{lemma}\label{abcohom2}
     Suppose that the class $[\cat{H}]$ vanishes in 
     $H^3(\pi_0(\cat{H});\on{A}_{\cat{H}})$. Then for any commutative group 
     $G$, there is a split exact sequence of abelian groups
     \begin{equation}\label{extgr2}
     \xymatrix{ 1 \ar[r] & {H^2(\pi_0(\cat{H});G)} \ar[r] & 
     {\pi_0(\catHom[\otimes](\cat{H},G[1]))} \ar[r] & 
     {\Hom[\pi_0(\cat{H})](\on{A}_{\cat{H}},G)}  \ar[r] & 1.}
     \end{equation}
 \end{lemma}

\begin{proof}
By definition of $\delta$, we clearly get the exact sequence \eqref{extgr2}. 
One possible way to show that it splits is the following. Since 
$[\cat{H}]$ is trivial in $H^3(\pi_0(\cat{H});\on{A}_{\cat{H}})$, the 
"essentially exact" sequence of $gr$-categories
$$
\xymatrix{ 1 \ar[r] & {\on{A}_{\cat{H}}[1]} \ar[r] & {\cat{H}} \ar[r] 
&  {\pi_0(\cat{H})[0]} \ar[r] & 1 }
$$ 
splits. This means that there is an equivalence of $gr$-categories
$\cat{H}\simeq \left( \on{A}_{\cat{H}}\to[1]\pi_0(\cat{H}) \right)$. 
Hence, a direct computation as in Lemma \ref{abcohom} shows that there is an 
isomorphism of groups
$$ 
\pi_0(\catHom[\otimes](\on{A}_{\cat{H}}\to[1]\pi_0(\cat{H}),G[1]))\simeq 
H^2(\pi_0(\cat{H});G)\times \Hom[\pi_0(\cat{H})](\on{A}_{\cat{H}},G).
$$
\end{proof}

For $\cat{H}=\Pi_1(\Omega X)$, we have $\pi_0(\Pi_1(\Omega X))=\pi_1(X)$ and
$\on{A}_{\Pi_1(\Omega X})=\pi_2(X)$ and the class $k_2(X)=[\Pi_1(\Omega X)]$ 
in $H^3(\pi_1(X);\pi_2(X))$ is the so-called first Postnikov $k$-invariant of 
$X$. Hence, using Lemma \ref{abcohom}, \ref{abcohom2} and Proposition 
\ref{H^2ab}, we get
\begin{cor}[Hopf's theorem for 2-cohomology]
     Let $X$ be connected and locally relatively 2-connected and $G$ a 
     commutative group.
     \begin{itemize}
     \item[(i)] There exists an exact sequence of abelian 
     groups
     $$
     \xymatrix{ 1 \ar[r] & {H^2(\pi_1(X);G)} \ar[r] & 
     {H^2(X;G_X)} \ar[r] &  
     {\Hom[\pi_1(X)](\pi_2(X),G)} \ar[r]^-{\delta} & {H^3(\pi_1(X);G)},}
     $$
     where $G$ is view as a $\pi_1(X)$-module with trivial action.
     
     \item[(ii)] If moreover the Postnikov invariant $k_2(X)$ vanishes in 
     $H^3(\pi_1(X);\pi_2(X))$, there is a split exact sequence of 
     abelian groups
     \begin{equation}\label{extab}
     \xymatrix{ 1 \ar[r] & {H^2(\pi_1(X);G)} \ar[r] & 
     {H^2(X;G_X)} \ar[r] & 
     {\Hom[\pi_1(X)](\pi_2(X),G)}  \ar[r] & 1.}
     \end{equation}
     \end{itemize} 
\end{cor}

\noindent
Recall that $H^2(\pi_1(X);G)\simeq \she_c (\pi_1(X);G)$, the group of 
equivalence classes of central extensions of $\pi_1(X)$ by $G$. Hence 
\eqref{extab} takes a similar form as the Universal Coefficient Theorem.

Now, let $G$ be a not necessarily commutative group and consider $\catc=G[1]$. 
Then the 2-category $\Cat{\Gamma}(X,\catEnd(G[1])[1]_X)$ is equivalent to 
$\Cat{Ger}(G_X)$, the strict 2-category of $G_X$-gerbes.
One may show (see for example \cite{Breen2}) that there is an 
isomorphism of pointed sets\footnote{The pointed set $\Two{\pi}_0(\Cat{Ger}(G_X))$ 
is sometimes denoted by $H_g^2 (X;G_X)$ and called the Giraud's second non 
abelian cohomology set of $X$ with values in $G_X$.}
$$
\Two{\pi}_0(\Cat{Ger}(G_X)) \simeq  H^1 (X;G_X\to[\ad]\Aut[\catGr](G_X)),
$$
where the right hand side is the first cohomology set of $X$ with values in 
the sheaf of crossed modules $G_X\to[\ad]\Aut[\catGr](G_X)$.

Assume that $X$ is locally relatively 2-connected. By Theorem 
\ref{mu2equiv}, there is an equivalence of 2-categories 
$$\Cat{Ger}(G_X) \overset{\sim}{\lra} \CatHom(\Pi_2(X),\catEnd(G[1])[1]).$$
Since $\Pic(G[1])\simeq \on{Out}_{\catGr}(G)$, from \eqref{mu2endCbis} we get
\begin{prop}[Hurewicz-Hopf's formula II]\label{H^2nonab}
     Let $X$ be connected and locally relatively 2-connected. 
     Then for any group $G$ there is an isomorphism of pointed sets 
     $$  
     H^1 (X;G_X\to[\ad]\Aut[\catGr](G_X))
     \isoto\pi_0(\catHom[\otimes](\Pi_1(\Omega X),G\to[\ad]\Aut[\catGr](G)))/
     \on{Out}_{\catGr}(G).
     $$
\end{prop}

With similar computations as for the commutative case, we get the following 
\begin{lemma}\label{nonabcohom}
     Let $\catg=G^{-1}\to[d]G^0$. Then \eqref{exgr} gives an exact sequence
     of pointed sets
     $$
     \xymatrix{ 1 \ar[r] & {H^1(\pi_0(\cat{H});G^{-1}\to[d]G^0)} \ar[r] & 
     {\pi_0(\catHom[\otimes](\cat{H},G^{-1}\to[d]G^0))} \ar[r] &  
     {\Hom[\catGr](\on{A}_{\cat{H}},\ker d)} ,}
     $$
     where the first term is the cohomolgy set of $\pi_0(\cat{H})$ with values 
     in the crossed module $G^{-1}\to[d]G^0$ (see for example \cite{Breen}).
\end{lemma} 
Combining Proposition \ref{H^2nonab} and Lemma \ref{nonabcohom} with 
$\cat{H} = \Pi_1(\Omega X)$ and $\catg= G\to[\ad]\Aut[\catGr](G)$, we get
\begin{cor}[Hopf's theorem for non abelian 2-cohomology]
     Let $X$ be connected and locally relatively 2-connected. 
     Then for each group $G$ there is an exact sequence of pointed sets
     $$
     \xymatrix{ 1 \ar[r] & {H^1(\pi_1(X);G\to[\ad]\Aut[\catGr](G))/
     \on{Out}_{\catGr}(G)} \ar[r] & 
     {H^1(X;G_X\to[\ad]\Aut[\catGr](G_X))} \lra&\\
     & & \lra {\Hom[\catGr](\pi_2(X),\on{Z}(G))/\on{Out}_{\catGr}(G),}  &}
     $$
     where the action of $\on{Out}_{\catGr}(G)$ on $H^1(\pi_1(X);G\to[\ad]
     \Aut[\catGr](G))$ is by conjugation, and that on $\Hom[\catGr](\pi_2 (X),
     \on{Z}(G))$ is induced by the natural action on $\on{Z}(G)$.
     If moreover $\pi_1(X)$ is trivial, one gets an isomorphism (``Hurewicz's 
     formula'')
     $$ H^1 (X;G_X\to[\ad]\Aut[\catGr](G_X)) \isoto
     \Hom[\catGr](\pi_2 (X),\on{Z}(G))/\on{Out}_{\catGr}(G),$$
\end{cor}
\noindent
A similar result holds replacing $G\to[\ad]\Aut[\catGr](G)$ by a general 
crossed module $G^{-1}\to[d]G^0$.

Recall that the pointed set $H^1(\pi_1(X);G\to[\ad]\Aut[\catGr](G))$ 
classifies extensions of $\pi_1(X)$ by $G$. Denoting by  $\she (\pi_1(X);G)$ 
the set of equivalence classes of such extensions and by $H_g^2 (X;G_X)$ the 
set $H^2(X;G_X\to[\ad]\Aut[\catGr](G_X))$, the previous sequence takes the 
more familiar form
$$
     \xymatrix{ 1 \ar[r] & {\she(\pi_1(X);G)/\on{Out}_{\catGr}(G)} \ar[r] & 
     {H_g^2(X;G_X)} \ar[r] & 
     {\Hom[\catGr](\pi_2(X),\on{Z}(G))/\on{Out}_{\catGr}(G)}.}
$$

\section*{Final comments}

What's next? It seems clear that, using the same technique, one should expect 
for each $n$-category $\Catc$ and each locally relatively $n$-connected space 
$X$ a natural $n$-equivalence
\begin{equation}\label{n-muequiv}
\Two{\mu^n}\colon \Cat{\Gamma}(X,\Catc _X)\overset{\sim}{\lra}
\CatHom[n\CatCat](\Pi_n(X),\Catc),
\end{equation}
where $\Cat{\Gamma}(X,\Catc _X)$ denotes the global sections of the constant 
$n$-stack with stalk $\Catc$, $\Pi_n(X)$ the homotopy $n$-groupoid of $X$ and 
$\CatHom[n\CatCat](\cdot,\cdot)$ the $n$-category of $n$-functors.
However, some care has to be taken since there are several non-equivalent 
definitions of $n$-categories for $n\geqs 3$ (see for example 
\cite{Leinster}). We will not investigate this problem any further here. 
An answer in this direction, using the formalism of Segal categories, is 
partially given in \cite{Toen} for $\Catc=(n-1)\CatCat$, the strict 
$n$-category of $(n-1)$-categories, and $X$ a pointed and connected 
$CW$-complex.

Let us suppose for a while that formula \eqref{n-muequiv} is valid and 
consider a commutative group $G$. Denote by $G[n]$ the strict 
$gr$-$n$-category with a single element, trivial $i$-arrows for $i\leq n-1$ 
and $G$ as $n$-arrows. Then one may check that there is an isomorphism of 
groups $$H^n (X;G_X)\simeq \Two{\pi}_0(\Cat{\Gamma}(X,G[n]_X)),$$
where the right-hand side is the group of $n$-equivalence classes of 
global objects in $G[n]_X$.\\
Suppose that $X$ is  locally relatively $n$-connected. From \eqref{n-muequiv}, 
we have an isomorphism of groups
\begin{equation}\label{n-Hopf}
H^n (X;G_X)\simeq\Two{\pi}_0(\CatHom[n\CatCat](\Pi_n(X),G[n])).
\end{equation}
\noindent
This isomorphism should be interpreted as the ``Hurewicz-Hopf's formula''. 
Indeed, if we suppose that $X$ is connected and that $\pi_i(X)\simeq 1$ for 
all $2 \leq  i \leq n-1$, we get an "essentially exact" sequence of 
$gr$-$(n-1)$-categories 
$$
\xymatrix{ 1 \ar[r] & {\pi_n(X)}[n-1] \ar[r] & {\Pi_{n-1}(\Omega X)} \ar[r] &  
{\pi_1(X)[0]} \ar[r] & 1 ,}
$$ 
and hence an exact sequence of groups
$$
\xymatrix{ 1 \ar[r] & {\Two{\pi}_0(\CatHom[\otimes](\pi_1(X)[0],G[n-1]))}  
     \ar[r] & {\Two{\pi}_0(\CatHom[\otimes](\Pi_{n-1}(\Omega X),G[n-1]))}
     \lra & \\
     & & \lra {\Two{\pi}_0(\CatHom[\otimes](\pi_n(X)[n-1],G[n-1]))} &}
$$
From the isomorphism \eqref{n-Hopf} and a direct calculation, we shall finally 
get the Hopf's exact sequence
$$
\xymatrix{ 1 \ar[r] & {H^n(\pi_1(X);G)} \ar[r] & {H^n(X;G_X)} \ar[r] & 
{\Hom[\catGr](\pi_n(X),G)} ,}
$$
where $G$ is view as a $\pi_1(X)$-module with trivial action
(we refer to \cite{Eilenberg-MacLane}, for a classical proof of this result).
If moreover $\pi_1(X)\simeq 1$, the Hurewicz's morphism 
$H^n(X;G_X)\to\Hom[\catGr](\pi_n(X),G)$ is an isomorphism.

If $G$ is a not necessarily commutative group, we define the 
$n$-groupoid $G[\![n]\!]$ by induction as
$$ 
G[\![1]\!]=G[1] , \qquad G[\![n+1]\!] = \Cat{Aut}_{n\CatCat}(G[\![n]\!])[1],
$$
where $\Cat{Aut}_{n\CatCat}(G[\![n]\!])$ denotes the $gr$-$n$-category 
of auto-$n$-equivalence of $G[\![n]\!]$ (note that, when $G$ is commutative, 
if we require $G$-linearity at each step in the definition of $G[\![n]\!]$, we 
recover $G[n]$). Then one may define the non abelian $n$-cohomology 
set of $X$ with coefficient in $G_X$ as 
$$H_g^n(X;G)=\Two{\pi}_0(\Cat{\Gamma}(X,G[\![n]\!]_X)).$$
If $X$ is locally relatively $n$-connected, then the $n$-equivalence 
\eqref{n-muequiv} gives an isomorphism of pointed sets    
$$H_g^n(X;G)\simeq\Two{\pi}_0(\CatHom[n\CatCat](\Pi_n(X),G[\![n]\!])).$$
This is the non abelian version of the ``Hurewicz-Hopf's formula''.

\begin{appendix}

\section{The stack of sheaves with values in a complete category}
\label{sheavesC}

We recall here the construction of the stack of sheaves with values in a 
complete category $\catc$, \emph{i.e.} a category which admits all small 
limits.

Let $X$ be a topological space and denote by $\catOp(X)$ the category of its 
open subset with inclusions morphisms. 
\begin{defi}
     A presheaf on $X$ with values in $\catc$ is a functor
     $$ \catOp(X)^{\op} \lra  \catc. $$ 
     A morphism between presheaves is a morphism of functors. We denote
     by $\catPSh_X(\catc)$ the category of presheaves on $X$ with values in 
     $\catc$.\\
     A presheaf is called a sheaf if it commutes to filtered limits indexed by
     coverings that are stable by finite intersection\footnote{Recall that, if 
     $\shf$ is a presheaf of sets the sheaf condition means that for any open 
     subset $U\subset X$, and any open covering $\{U_i\}_{i\in I}$ of $U$, 
     the natural sequence given by the restriction maps
     $$
     \xymatrix{ \shf(U) \ar[r] & \prod_{i\in I} \shf(U_i) \ar@<.2pc>[r]
     \ar@<-.2pc>[r] & \prod_{i,j\in I} \shf(U_{ij}) }
     $$
     is exact in the usual sense.}, and we denote by $\catSh_X(\catc)$ the 
     full subcategory of $\catPSh_X(\catc)$ whose objects are sheaves.
\end{defi}

Note that if $U\subset X$ is an open subset and $\shf$ is a sheaf on $X$, then
its restriction $\shf|_{U}$ is also a sheaf. Hence we can define the prestack 
of sheaves on $X$, denoted by $\stkSh_X(\catc)$, by assigning 
$X\supset U\mapsto \catSh_U(\catc)$.

Let $\pshf,\pshg$ be two presheaves on $X$. We have a natural 
bijective map of sets
$$ 
\Hom[\catPSh_X(\catc)](\pshf,\pshg)\isoto
\underset{(U,V) \atop V\subset U}{\varprojlim}{\Hom[\catc](\pshf(U),\pshg(V))} 
$$
where $(U,V)$ is considered as an object of $\catOp(X)^{\op}\times
\catOp(X)$. Now let $U\subset X$ be an open subset, $\pshf$ a presheaf on 
$X$ and $\shg$ a presheaf on $U$. Then it is easy to see that we have the 
isomorphism of sets
$$ 
\Hom[\catPSh_U(\catc)](\pshf|_{U},\pshg)\isoto
\underset{(V,W) \atop W\subset V\subset U}{\varprojlim}
{\Hom[\catc](\pshf(V),\pshg(W))} \isoto
\underset{(V,W) \atop W\subset V\subset X}{\varprojlim}
{\Hom[\catc](\pshf(V),\pshg(W\cap U))}.
 $$

\begin{lemma}\label{homlemma}
     Let $\pshf$ be a presheaf and $\shg$ a sheaf on $X$.
     Then the presheaf $\ \ $ $\shhom[\catPSh_X(\catc)](\pshf,\shg)$ defined 
     by 
     $$ 
     X\supset U\mapsto \shhom[\catPSh_X(\catc)](\pshf,\shg)(U)=
     \Hom[\catPSh_U(\catc)](\pshf|_U,\shg|_U) 
     $$
     is a sheaf of sets.
\end{lemma}
 \begin{proof}
We have to show that $\shhom[\catPSh_X(\catc)](\pshf,\shg)$ commutes to small 
filtered limits indexed by coverings that are stable by finite intersection. 
Let $\{U_i\}_{i\in I}$ be such a covering of an open subset $U\subset X$.
Then we have
\begin{align*}
   \shhom[\catPSh_X(\catc)](\pshf,\shg)(U) & =
   \Hom[\catPSh_U(\catc)](\pshf|_U,\shg|_U)
   \simeq \underset{(V,W) \atop W\subset V}{\varprojlim}
   \Hom[\catc](\pshf(V),\shg(U\cap W))\\
   & \simeq   \underset{(V,W) \atop W\subset V}{\varprojlim}
    \Hom[\catc](\pshf(V),\underset{i\in I}{\varprojlim}\shg(U_i\cap W))\\
   & \simeq  \underset{(V,W) \atop W\subset V}{\varprojlim}
   \underset{i\in I}{\varprojlim}
   \Hom[\catc](\pshf(V),\shg(U_i\cap W))\\
   & \simeq  \underset{i\in I}{\varprojlim}
   \underset{(V,W) \atop W\subset V}{\varprojlim}
   \Hom[\catc](\pshf(V),\shg(U_i\cap W))\\
   &\simeq  \underset{i\in I}{\varprojlim}
   \Hom[\catPSh_{U_i}(\catc)](\pshf|_{U_i},\shg|_{U_i})
   \simeq \underset{i\in I}{\varprojlim}
   \shhom[\catPSh_X(\catc)](\pshf,\shg)(U_i).
\end{align*} 
\end{proof}

\begin{lemma}\label{techlemma0}
     Let $\pshf$ be a presheaf on $X$. Then $\pshf$ is a sheaf if and only if 
     for any object $A\in\Ob{\catc}$ and any open subset $U\subset X$ the 
     presheaf
     $$ U\supset V\mapsto \Hom[\catc](A,\pshf(V)) $$
     is a sheaf of sets.
\end{lemma}

\begin{proof}
Follows immediately from Yoneda's Lemma.
\end{proof}

\begin{prop}
     The prestack $\stkSh_X(\catc)$ of sheaves with values in $\catc$ is a 
     stack. 
\end{prop}

\begin{proof}
By Lemma \ref{homlemma}, the prestack is seperated. Now let 
$(\{U_{i}\}_{i\in I}, \{\shf_{i}\}_{i\in I}, \{\theta_{ij}\}_{i,j\in I})$
be a descent datum for $\stkSh_X(\catc)$ on open subset $U\subset X$. 
By taking a refinement, we can assume that the covering $\{U_{i}\}_{i\in I}$ 
is stable by finite intersections.\\
Let $V\subset U$. Then the cocycle condition allows us to define
$$ \shf(V)=\underset{i\in I}{\varprojlim}\shf_i(V\cap U_i). $$
It is then obvious that $\shf$ is a sheaf (for instance using Lemma 
\ref{techlemma0} and the fact that this is true if $\catc=\catSet$) 
which by construction is isomorphic to $\shf_i$ on $U_i$ by an isomorphism
$\theta_i$ such that $\theta_{ij}\circ\theta_j=\theta_i$ on $U_i\cap U_j$.
\end{proof}

\begin{prop}
    The stack $\stkSh_X(\catc)$ admits all small limits. 
\end{prop}

\begin{proof}
Let $\beta\colon\cati \ra \catSh_X(\catc)$ be a functor, with $\cati$ a small 
category. Then, for each open subset $U\subset X$, set
$$ \shf(U)=\underset{i\in \cati}{\varprojlim}{\beta(i)(U)}. $$
It is immediately verified that $\shf$ is a sheaf on $X$ that satisfies
$\shf \simeq \underset{i\in \cati}{\varprojlim}{\beta(i)}$.
\end{proof}

\begin{defi}
     Let $\pshf$ be a presheaf. A sheaf $\wtilde{\pshf}$ together with a 
     morphism $\pshf \ra \wtilde{\pshf}$ is called the sheaf associated to 
     $\pshf$ if it satisfies the usual universal property, \emph{i.e.} any 
     morphism from $\pshf$ into a sheaf $\shg$ factors uniquely through 
     $\wtilde{\pshf}$:
     $$ \xymatrix{ 
     {\pshf}  \ar[r] \ar[d] & {\shg} \\
     {\wtilde{\pshf}} \ar[ur] & } 
     $$
\end{defi}

\begin{prop}\label{const}
     Let $M\in\Ob{\catc}$. Assume that $X$ is locally connected. 
     Then the sheaf associated to the constant presheaf with stalk $M$ exists.
\end{prop}

\begin{proof}
Let $U\subset X$ be an open subset. Denote by $\#U$ the set of connected 
components of $U$. Set
$$ M_X(U)=M^{\#U}. $$
Let $x\in U$. Then we denote by $\ol{x}_U$ its class in $\#U$.\\
For any inclusion $V\subset U$ of open subsets and for any $x\in V$,
we have the natural morphisms 
$$ M_X(U) \ra M_{\ol{x}_U}\ra M_{\ol{x}_V}. $$
These morphisms define the morphism
$$ M_X(U) \ra M_X(V). $$
Since we know that this is a sheaf if $\catc=\catSet$, Lemma \ref{techlemma0} 
implies that $M_X$ is a sheaf which verifies the desired universal property.
\end{proof}

\begin{defi}
     Let $M\in\Ob{\catc}$. The sheaf associated to the constant presheaf with 
     stalk $M$ is called the constant sheaf with stalk $M$, and we denote
     it by $M_X$.
\end{defi}

\begin{remark}
The hypothesis on the connectivity of $X$ in the Proposition \ref{const} is 
necessary to recover the classical definition of constant sheaf. 
More precisely, if $M$ is a set and $M_X$ is the constant sheaf 
defined in the usual way, there is a natural injective map (the 0-monodromy)
$\mu^0\colon M_X(X)\to \Hom(\#X,M)\simeq M^{\#X}$ defined by $\mu^0(s)(\ol{x}_X) = 
s(\ol{x}_X)$ for a section $s$ in $M_X(X)$. Clearly, if $X$ is locally 
connected, $\mu^0$ is a bijection.
\end{remark}

Let $X$ be a locally connected topological space. 
Denote by $\catCSh_X(\catc)$ the full subcategory of $\catSh_X(\catc)$ of 
constant sheaves. The previous construction defines a faithful functor
$$(\cdot)_X\colon \catc \lra \catCSh_X(\catc) ,$$
which is an equivalence if $X$ is connected (a quasi-inverse is given by the 
global sections functor).

\begin{defi}
     A sheaf $\shf$ is called locally constant if there is an open covering
     $X=\bigcup {U_i}$ such that $\shf|_{U_i}$ is isomorphic to a constant 
     sheaf.\\
     We denote by $\stkLoc_X(\catc)$ the full substack of $\stkSh_X(\catc)$ 
     whose objects are the locally constant sheaves.
\end{defi}

\section{The 2-stack of stacks with values in a 2-complete 2-category}
\label{stacksC}

Let us recall the construction of the 2-stack of stacks with values in a 
2-complete 2-category $\Catc$, \emph{i.e.} a 2-category which admits all small 
2-limits (references for the basic definitions about 2-stacks are made to 
\cite{Breen}).
Recall that the 2-Yoneda lemma states that the 2-functor
$$ 
\Catc \lra \widehat\Catc = \CatHom(\Catc^{\op},\CatCat) \quad ; 
\quad \on{P} \mapsto \catHom[\Catc](\,\cdot\,,\on{P}) 
$$
is fully faithful (see for example \cite[cap. 1]{Leinster}, for more details). 
Since $\widehat\Catc$ is strict, the reader may assume for sake of simplicity 
that $\Catc$ is a strict 2-category.

Let $X$ be a topological space and denote by $\Cat{Op}(X)$ the 2-category of
its open subsets, obtained by trivially enriching $\cat{Op}(X)$ with identity 
2-arrows.

\begin{defi}
     A prestack on $X$ with values in $\Catc$ is a 2-functor
     $$ \Cat{Op}(X)^{\op} \lra  \Catc. $$ 
     A functor between prestacks is a 2-transformation of 2-functors and 
     transformations of functors of prestacks are modifications of 
     2-transformations of 2-functors. We denote by $\CatPSt_X(\Catc)$ 
     the 2-category of prestacks on $X$ with values in $\Catc$.\\
     A prestack is called a stack if it commutes to filtered 2-limits indexed 
     by coverings that are stable by finite intersection\footnote{Similarly to 
     the definition of sheaf, if $\stks$ is a prestack of categories the 
     stack condition means that for any open subset $U\subset X$, and any open 
     covering $\{U_i\}_{i\in I}$ of $U$, the natural sequence given by the 
     restriction functors
     $$
     \xymatrix{ \stks(U) \ar[r] & \prod_{i\in I} \stks(U_i) \ar@<.2pc>[r]
     \ar@<-.2pc>[r] & \prod_{i,j\in I} \stks(U_{ij}) \ar[r] \ar@<.3pc>[r]
     \ar@<-.3pc>[r] & \prod_{i,j,k\in I} \stks(U_{ijk}) }
     $$
     is exact in the sense of~\cite[expos\'e XIII]{SGA1}.}, and we denote by 
     $\CatSt_X(\Catc)$ the full subcategory of $\CatPSt_X(\Catc)$ whose 
     objects are stacks.
\end{defi}

Note that if $U\subset X$ is an open subset and $\stks$ is a stack on $X$, then
its restriction $\stks|_{U}$ is also a stack. Hence the assignment 
$X\supset U\mapsto \CatSt_U(\Catc)$ defines the pre-2-stack of stacks on $X$ 
with values in $\Catc$, which we denote by $\twostkSt_X(\Catc)$.\\
Let $\pstks,\pstkt$ be two prestacks on $X$. We have a natural 
equivalence of categories
$$ 
\catHom[\CatPSt_X(\Catc)](\pstks,\pstkt)\overset{\sim}{\lra}
\twolim {(U,V) \atop V\subset U} {\catHom[\Catc](\pstks(U),\pstkt(V))} 
$$
where $(U,V)$ is considered as an object of $\Cat{Op}(X)^{\op}\times
\Cat{Op}(X)$. Now let $U\subset X$ be an open subset, $\pstks$ a prestack on 
$X$ and $\pstkt$ a prestack on $U$. Then it is easy to see that we have the 
equivalence of categories
$$ 
\catHom[\CatPSt_U(\Catc)](\pstks|_{U},\pstkt)\overset{\sim}{\lra}
\twolim {(V,W) \atop W\subset V\subset U}
{\catHom[\Catc](\pstks(V),\pstkt(W))} \overset{\sim}{\lra}
\twolim{(V,W) \atop W\subset V\subset X}
{\catHom[\Catc](\pstks(V),\pstkt(W\cap U))}.
 $$

Similarly to the case of sheaves, we have 
\begin{lemma}\label{homlemmab}
     Let $\pstks$ be a prestack and $\stkt$ be a stack on $X$.\\
     Then the prestack $\stkHom[\CatPSt_X(\Catc)](\pstks,\stkt)$ defined by 
     $$ 
     \stkHom[\CatPSt_X(\Catc)](\pstks,\stkt)(U)=
     \catHom[\CatPSt_U(\Catc)](\pstks|_U,\stkt|_U) 
     $$
     is a stack of categories.
\end{lemma}


Moreover, using the 2-Yoneda's Lemma, we get
\begin{lemma}\label{techlemmab0}
     Let $\pstks$ be a prestack on $X$. Then $\pstks$ is a stack if and only if
     for any object $\on{P}\in\Ob{\Catc}$ and any open subset $U\subset X$ the 
     prestack
     $$ U\supset V\mapsto \catHom[\Catc](\on{P},\pstks(V)) $$
     is a stack of categories.
\end{lemma}

\begin{prop}
     The pre-2-stack $\twostkSt_X(\Catc)$ of stacks with values in $\Catc$ is 
     a 2-stack. 
\end{prop}

\begin{proof}
By Lemma \ref{homlemmab}, the pre-2-stack is seperated. Now let 
$$
(\{U_{i}\}_{i\in I}, \{\stks_{i}\}_{i\in I}, \{F_{ij}\}_{i,j\in I}, 
\{\varphi_{ijk}\}_{i,j,k\in I})
$$
be a descent datum for $\twostkSt_X(\Catc)$ on open subset $U\subset X$. 
This means that $\{U_i\}_{i\in I}$ is an open covering of $U$, $\stks_i$ are
stacks on $U_i$, $F_{ij}\colon \stks_j\vert_{U_{ij}}\overset{\sim}{\lra}
\stks_i\vert_{U_{ij}}$ are equivalences of stacks and 
$\varphi_{ijk}\colon \varphi_{ij}\circ \varphi_{jk} \to \varphi_{ik}$ are 
invertible transformations of functors from $\stks_k\vert_{U_{ijk}}$ to 
$\stks_i\vert_{U_{ijk}}$, such that for any $i,j,k,l \in I$, the following 
diagram of transformations of functors from $\stks_l|_{U_{ijkl}}$ to 
$\stks_i|_{U_{ijkl}}$ commutes
\begin{equation}\label{gluestk}
\xymatrix@C=4em{ {F_{ij}\circ F_{jk}\circ F_{kl}}
\ar[r]^-{\varphi_{ijk}} \ar[d]^{\varphi_{jkl}}
&{F_{ik}\circ F_{kl}}\ar[d]^{\varphi_{ikl}}\\
{F_{ij}\circ F_{jl}}\ar[r]^{\varphi_{ijl}}&F_{il}. }
\end{equation}
By taking a refinement, we can assume that the covering $\{U_{i}\}_{i\in I}$ 
is stable by finite intersections.\\
Let $V\subset U$ be an open subset. Then the cocycle condition \eqref{gluestk} 
allows us to define the category
$$ \stks(V)=\twolim {i\in I}{\stks_i(V\cap U_i)}. $$
It is then obvious that the assignment $U\supset V\mapsto \stks (V)$ defines 
a stack $\stks$ (for instance using Lemma \ref{techlemmab0} and the fact that 
this is true if $\Catc=\CatCat$) and that, by construction, there are e
quivalences of stacks $F_i\colon \stks\vert_{U_i}\overset{\sim}{\lra}\stks_i$. 
Moreover, one checks that there exist invertible transformations of functors 
$\varphi_{ij} \colon F_{ij}\circ F_j\vert_{U_{ij}} \isoto F_i\vert_{U_{ij}}$  
such that $\varphi_{ij}|_{U_{ijk}} \circ \varphi_{jk}|_{U_{ijk}} = 
\varphi_{ik}|_{U_{ijk}} \circ \varphi_{ijk}$.  
\end{proof}

\begin{prop}
    The 2-stack $\twostkSt_X(\Catc)$ admits all small 2-limits. 
\end{prop}

\begin{proof}
Let $\beta\colon\Cat{I} \ra \CatSt_X(\Catc)$ be a 2-functor, with $\Cat{I}$ 
a small 2-category. Then, for each open subset $U\subset X$, set
$$ \stks(U)=\twolim{i\in \Cat{I}}{\beta(i)(U)}. $$
It is immediately verified that $\stks$ is a stack on $X$ that satisfies
$\stks \simeq \twolim{i\in \Cat{I}}{\beta(i)}$.
\end{proof}

\begin{defi}
     Let $\pstks$ be a prestack. A stack $\wtilde{\pstks}$ together with a 
     functor $\pstks \lra \wtilde{\pstks}$ is called the stack associated to 
     $\pstks$ if it satisfies the usual universal property, \emph{i.e.} any 
     functor from $\pstks$ into a stack $\stkt$ factors through 
     $\wtilde{\pstks}$ up to unique equivalence:
     $$ \xymatrix{ 
     {\pstks}  \ar[r] \ar[d] & {\stkt} \\
     {\wtilde{\pstks}} \ar[ur] & } 
     $$
\end{defi}

\begin{prop}
     Let $\on{P}\in\Ob{\Catc}$. Assume that $X$ is locally 1-connected. 
     Then the stack associated to the constant prestack with stalk $\on{P}$ 
     exists.
\end{prop}

\begin{proof}
Let $U\subset X$ be an open subset. Set
$$ \on{P}_X(U)=\on{P}^{\Pi_1(U)}, $$
where $\on{P}^{\Pi_1(U)}$ denotes the 2-limit of the constant 2-functor
$\cat{\Delta}(P)\colon \Pi_1(U)\lra \Catc$ at $\on{P}$.
Let $x\in V\subset U$ and denote by $x_U$ the image of $x$ in $\Pi_1(U)$ by the
natural functor $\Pi_1(V)\lra\Pi_1(U)$. Set $\on{P}_{x_U}=\cat{\Delta}(P)(x_U)$
and similarly for $\on{P}_{x_V}$. Then we have the natural 1-arrows in $\Catc$
$$ \on{P}_X(U) \lra \on{P}_{x_U}\lra \on{P}_{x_V},$$
which define the 1-arrow
$$ \on{P}_X(U) \lra \on{P}_X(V). $$
If $\Catc=\CatCat$, then $\on{P}$ is a category and by Theorem \ref{muequiv} 
there are equivalences $\on{P}^{\Pi_1(U)}\simeq\catHom (\Pi_1(U),\on{P})
\underset{\mu}{\overset{\sim}{\longleftarrow}}\on{P}_U(U)$, hence the 
assignment $X\supset U\mapsto \on{P}^{\Pi_1(U)}$ defines the stack of locally 
constant sheaves on $U$ with values in $\on{P}$. In the general case, one use 
Lemma \ref{techlemmab0} to show that this construction gives a stack which 
verifies the desired universal property.
\end{proof}

\begin{defi}
     Let $\on{P}\in\Ob{\Catc}$. The stack associated to the constant prestack 
     with stalk $\on{P}$ is called the constant stack with stalk $\on{P}$, and 
     we denote it by $\on{P}_X$.
\end{defi}

Let $X$ be a locally 1-connected topological space. 
Denote by $\CatCSt_X(\Catc)$ the full sub-2-category of $\CatSt_X(\Catc)$ of 
constant stacks. The previous construction defines a faithful 2-functor
$$(\cdot)_X\colon \Catc \lra \CatCSt_X(\Catc) ,$$
which is an equivalence if $X$ is 1-connected (a quasi-2-inverse is given by 
the global sections 2-functor). This easily follows  from Corollary 
\ref{cstacks}.

\begin{defi}
     A stack $\stks$ is called locally constant if there exists an open 
     covering $X=\bigcup {U_i}$ such that $\stks|_{U_i}$ is isomorphic to a 
     constant stack.\\
     We denote by $\twostkLoc_X(\Catc)$ the full sub-2-stack of 
     $\twostkSt_X(\Catc)$ whose objects are the locally constant stacks.
\end{defi}

\end{appendix}

\providecommand{\bysame}{\leavevmode\hbox to3em{\hrulefill}\thinspace}

\vspace*{1cm}

\parbox[t]{15em}
{\scriptsize{
\noindent
Pietro Polesello\\
Universit{\`a} di Padova\\
Dipartimento di Matematica\\
via G. Belzoni, 7
35131 Padova, Italy\\
or: Universit{\'e} Pierre et Marie Curie\\
Institut de Math{\'e}matiques\\
175, rue du Chevaleret, 75013 Paris France\\
pietro@math.jussieu.fr\\}}
\qquad
\parbox[t]{15em}
{\scriptsize{
Ingo Waschkies\\
Laboratoire J.A. Dieudonn\'e\\
Universit\'e de Nice Sophia-Antipolis\\
Parc Valrose, 06108 Nice Cedex 2\\
ingo@math.unice.fr\\}}

\end{document}